\documentclass{amsart}
\usepackage[english]{babel}
\usepackage[utf8x]{inputenc}
\usepackage[T1]{fontenc}
\usepackage{amssymb}
\usepackage{bigints}
\numberwithin{equation}{section}
\usepackage{amsmath}
\usepackage{graphicx}
\usepackage[colorinlistoftodos]{todonotes}
\usepackage{caption}
\usepackage{subcaption}
\usepackage{bm}% bold math
\usepackage{amsmath,amsthm,amssymb,stmaryrd}
\usepackage{enumerate}

\newcommand{\bbint}[2]{\ensuremath{\;\backslash\!\!\!\!\backslash\!\!\!\!\!\int_{#1}^{#2}}}

\newmuskip\pFqmuskip
\newcommand*\pFq[6][8]{%
	\begingroup % only local assignments
	\pFqmuskip=#1mu\relax
	\mathchardef\normalcomma=\mathcode`,
	\mathcode`\,=\string"8000
	\begingroup\lccode`\~=`\,
	\lowercase{\endgroup\let~}\pFqcomma
	{}_{#2}F_{#3}{\left(\genfrac..{0pt}{}{#4}{#5}\Bigg| #6\right)}%
	\endgroup
}
\newcommand{\pFqcomma}{{\normalcomma}\mskip\pFqmuskip}

\title[Finite-Part Integration]{Finite-Part Integration in the Presence of Competing Singularities: Transformation Equations for the hypergeometric functions arising from Finite-Part Integration}
\author{Lloyd Villanueva and Eric A. Galapon}

\address{Theoretical Physics Group, National Institute of Physics\\University of the Philippines, Diliman Quezon City\\1101 Philippines}

\email{eagalapon@up.edu.ph}
\date{\today}

\begin{document}
\maketitle

\begin{abstract}
Finite-part integration is a recently introduced method of evaluating convergent integrals by means of the finite part of  divergent integrals [E.A. Galapon, {\it Proc. R. Soc. A 473, 20160567} (2017)]. Current application of the method involves exact and asymptotic evaluation of the generalized Stieltjes transform $\int_0^a f(x)/(\omega + x)^{\rho} \, \mathrm{d}x$ under the assumption that the extension of $f(x)$ in the complex plane is entire. In this paper, the method is elaborated further and extended to accommodate the presence of competing singularities of the complex extension of $f(x)$. Finite part integration is then applied to derive consequences of known Stieltjes integral representations of the Gauss function and the generalized hypergeometric function which involve Stieltjes transforms of functions with complex extensions having singularities in the complex plane. Transformation equations for the Gauss function are obtained from which known transformation equations are shown to follow. Also, building on the results for the Gauss function, transformation equations involving the generalized hypergeometric function $\,_3F_2$ are derived.
\end{abstract}

\section{Introduction}
Divergent integrals arise in many areas of applied mathematics, such as in physics \cite{collins} and in engineering \cite{ang}. In those applications, the problem typically is how to assign meaningful physical values to them or compute the values of particular assignments to the divergent integrals, such as their finite parts. Recently, in revisiting the problem of missing terms arising from term by term integration leading to an infinite series of divergent integrals \cite{wong,distribution,mcwong}, it was determined that the finite part of divergent integrals can be rigorously used as a means of evaluating convergent integrals, a method we have referred to as finite-part integration \cite{galapon2}. Finite-part integration has been applied in the exact and asymptotic evaluation of the Stieltjes transform of integer \cite{tica1} and non-integer orders \cite{tica2}. Applying the method to known Stieltjes integral representations of some special functions has led to new representations of them. Moreover, more refined asymptotics of the special functions have been obtained from their new representations.  

The essential feature of finite-part integration is to give meaning to term by term integration that leads to divergent integrals. Consider the 
generalized Stieltjes transform \cite{widder,joshi,schwarz}
\begin{equation}\label{stielyjes}
    \mathcal{F}(\omega) =\int_0^a \frac{f(x)}{(\omega+x)^\rho}\, \mathrm{d}x, 
\end{equation}
where the upper limit of integration $a$ takes on the possible values $0<a\leq \infty$ and $\rho$ is a positive number, $\omega$ is generally complex with $|\mathrm{arg}\,\omega|<\pi$, and $f(x)$ is some integrable function that may depend on several parameters. One may attempt at evaluating this integral by expanding the kernel $(\omega+x)^{-\rho}$ binomially about $\omega=0$,
\begin{equation}\label{expand}
    \frac{1}{(\omega+x)^{\rho}}=\sum_{k=0}^{\infty} {-\rho\choose k} \frac{\omega^k}{x^{k+\rho}},
\end{equation}
inside the integral, and then distributing the integration over the summation to yield the infinite series of integrals,
\begin{equation}\label{series}
    \sum_{k=0}^{\infty} {-\rho \choose k} \omega^{k} \int_{0}^{a} \frac{f(x)}{x^{k+\rho}} \, \mathrm{d}x.
\end{equation}
If $f(x)$ has a zero of finite order at the origin, then the infinite series will eventually degenerate into a series of divergent integrals. The appearance of divergent integrals simply means that the interchange of summation and the integration are not allowed. Finite-part integration makes sense of this interchange, yielding the result  
\begin{equation}\label{byfpi}
    \int_0^a \frac{f(x)}{(\omega+x)^\rho}\, \mathrm{d}x = \sum_{k=0}^{\infty} {{-\rho}\choose{k}} \omega^k \bbint{0}{a} \frac{f(x)}{x^{k+\rho}}\,\mathrm{d}x + \Delta(\omega) ,
\end{equation}
where the integral $\bbint{0}{a}f(x)x^{-k-\rho}\,\mathrm{d}x$ is the finite part \cite{hadamard,monegato,galapon1} of the divergent integral $\int_0^a f(x)x^{-k-\rho}\,\mathrm{d}x$ and the term $\Delta(\omega)$ is a contribution coming from the singularity of the kernel of transformation at $-\omega$.
However, in \cite{galapon2, tica1, tica2} it was assumed that $f(x)$ possesses an entire complex extension $f(z)$. This condition severely restricts the domain of applicability of the method. 

In this paper, we elaborate the method of finite-part integration further to accommodate the presence of poles and branch points of the complex extension of $f(x)$. 
To that end, we will investigate the consequences of known integral representations of the Gauss function $_2F_1$ and the generalized hypergeometric function $_pF_q$ by direct finite part integration of the integrals in the representations. For the Gauss function, we will consider the known integral representation
\begin{equation}\label{mainlemma}
\pFq{2}{1}{\mu,\nu}{\mu+\rho}{1 - \frac{b}{a}} = \frac{\Gamma(\mu+\rho)}{\Gamma(\nu) \Gamma(\mu-\nu+\rho)} \frac{a^{\mu}}{b^{\nu-\rho}}\,  \int_0^{\infty} \frac{x^{\nu-1} }{(a+x)^{\mu}(b+x)^{\rho}} \, \mathrm{d}x,
\end{equation}
for $\mathrm{Re} (\nu)>0$, $\mathrm{Re} (\rho+\mu-\nu)>0$, $|\mathrm{arg} \, a|<\pi$ and $|\mathrm{Arg} \, b|<\pi$ \cite{saxena}. The integral in the right hand side can be interpreted as a Stieltjes transform of the function $f(x)=x^{\nu-1} (a+x)^{-\mu}$ with the kernel $(b+x)^{-\rho}$. The function $f(x)$ for this case has a complex extension that can have branch point at the origin, and a branch point or pole singularity at $-a$. For the generalized hypergeometric function, we will consider the integral representation
\begin{equation}\label{lukeq}
\, _{p+1}F_{q+1}\!\left.\left(\begin{array}{c}
\beta,\alpha_p\\
\beta+\alpha,\rho_q
\end{array}\right|z\right) = \frac{\Gamma(\beta+\sigma)}{\Gamma(\beta)\Gamma(\sigma)} \int_0^1 t^{\beta-1} (1-t)^{\sigma-1} \, _{p}F_{q}\!\left.\left(\begin{array}{c}
\alpha_p\\
\rho_q
\end{array}\right|z t\right)\,\mathrm{d}t  
\end{equation}
which is valid for $p\leq q+1$, $\mathrm{Re}(\beta)>0$, $\mathrm{Re}(\sigma)>0$, $|z|<1$ if $p=q+1$ \cite[p. 58, eqn. 10]{luke}. The integral can be cast as a Stieltjes integral with the substitution 
$t=s/(1+s)$,  
\begin{equation}\label{general}
\, _{p+1}F_{q+1}\!\left.\left(\begin{array}{c}
\beta,\alpha_p\\
\beta+\alpha,\rho_q
\end{array}\right|z\right) = \frac{\Gamma(\beta+\sigma)}{\Gamma(\beta)\Gamma(\sigma)} \int_0^{\infty} \frac{s^{\beta-1}}{(s+1)^{\beta+\sigma}} \, _{p}F_{q}\!\left.\left(\begin{array}{c}
\alpha_p\\
\rho_q
\end{array}\right|\frac{z s}{s+1}\right)\mathrm{d}s 
\end{equation}
under the same conditions on the parameters. This representation can be interpreted as a Stieltjes transform of the function 
\begin{equation*}
f(s)=s^{\beta-1}\, _{p}F_{q}\!\left.\left(\begin{array}{c}
\alpha_p\\
\rho_q
\end{array}\right|\frac{z s}{s+1}\right),
\end{equation*}
with the kernel $(s+t)^{-\beta-\sigma}$, evaluated at the specific point $t=1$. The complex extension of $f(s)$ may have a branch point at the origin, and a branch point or a pole at $-(1-z)^{-1}$. 

The Stieltejs integral representations \eqref{mainlemma} and \eqref{general} involve the transform of functions that have singularities competing with the singularity of the kernel of transformation, a circumstance which is outside the scope of \cite{galapon2,tica1,tica2}. Here we will perform finite-part integration in the presence of such singularities. For the Gauss function, we will show that finite-part integration of \eqref{mainlemma} under different values of the parameters leads to transformation equations in argument of $(1-z)$, and obtain from them known transformation equations of the Gauss function. For the generalized hypergeometric function, we will apply finite part integration on \eqref{general} for the specific case of the hypergeometric function $_3F_2$ for a large class of the parameters, and arrive at new transformation equations involving $_3F_2$.

The rest of the paper is organized as follows. In Section-\ref{finitepartintegration} we summarize the method of finite part integration. In Section-\ref{fundamental}, we obtain the pair of fundamental finite parts integrals arising from a branch point and pole singularity at the origin. In Section-\ref{origin}, we perform finite part integration on the Stieltjes integration representation \eqref{mainlemma} and obtain transformation equations in powers of $(1-z)$. In Section-\ref{section3f2}, we perform finite part integration on the representation \eqref{general} for the generalized hypergeometric function $_3F_2$ and obtain transformation equations in powers of $(1-z)$.  In Section-\ref{conclusion}, we conclude.  

Throughout the paper, $\mathbb{Z}$ denotes the set of integers; $\mathbb{Z}^+$, the positive integers; $\mathbb{Z}^-$, the negative integers; $\mathbb{Z}^+_0$, the positive integers including $0$; $\mathbb{Z}^-_0$, the negative integers including $0$.

\section{Finite Part Integration}\label{finitepartintegration}
Finite part integration is a method of evaluating convergent integrals by means of the finite part of divergent integrals. The method proceeds in two major steps. The first step is deliberately inducing divergent integrals from a given convergent integral, and the second step is recasting the given integral into a form that leads to its evaluation in terms of the finite parts of the induced divergent integrals. In general, other divergent integrals could arise in the process; for this reason, we refer to the divergent integrals that arise in the first step as fundamental and those that follow as progenic (a child of the fundamental divergent integral). The first step, once the fundamental divergent integrals have been identified, involves extracting the finite part of the divergent integrals and representing the finite parts as contour integrals in the complex plane. The second step involves extracting the given convergent integral from a contour integral whose characteristics are dictated by the nature of the fundamental divergent integrals. The transition from the first to the second step is made possible by the complex contour integral representation of the fundamental finite part integrals. We will detail in a short while how finite part integration will be implemented for our present problem. 

For the finite-part integration of the Stieltjes transform \eqref{stielyjes}, the fundamental divergent integrals are induced by expanding the kernel of transformation that leads to the infinite series of divergent integrals \eqref{series}. Then the fundamental divergent integrals are those of the form
\begin{equation}\label{divinteg}
    \int_0^a \frac{f(x)}{x^{\lambda}}\,\mathrm{d}x, \;\;\; \lambda\geq 1, \; 0<a\leq\infty,
\end{equation}
where the divergence arises only from a non-integrable singularity at the origin. For a finite upper limit of integration $a$, the finite part is obtained by temporarily removing the divergence by replacing the offending lower limit of integration with some positive $\epsilon<a$ and then grouping the resulting integral into a pair of terms,
\begin{equation}
    \int_{\epsilon}^a \frac{f(x)}{x^{\lambda}}\, \mathrm{d}x = C_{\epsilon} + D_{\epsilon},
\end{equation}
where $C_{\epsilon}$ is the group of terms that converges in the limit as $\epsilon\rightarrow 0$ and $D_{\epsilon}$ is the group of terms that diverge in the same limit. The finite part is obtained by dropping the diverging part $D_{\epsilon}$, leaving only the converging $C_{\epsilon}$ and assigning its limit as the value of the divergent integral. That is the finite part is given by
\begin{equation}\label{finitepart}
\bbint{0}{a} \frac{f(x)}{x^{\lambda}}\mathrm{d} x = \lim_{\epsilon\rightarrow 0} C_{\epsilon} .
\end{equation}
Equivalently the finite part also assumes the representation
\begin{equation}\label{finitepart2}
\bbint{0}{a} \frac{f(x)}{x^{\lambda}}\mathrm{d} x = \lim_{\epsilon\rightarrow 0} \left(\int_{\epsilon}^a \frac{f(x)}{x^{\lambda}}\,\mathrm{d}x-D_{\epsilon}\right) .
\end{equation}
By definition, the limits \eqref{finitepart} and \eqref{finitepart2} always exist and their value may be equal to zero under certain circumstances. When the upper limit of integration extends to infinity, $a=\infty$, the finite part is given by the limit
\begin{equation}\label{finitepart3}
\bbint{0}{\infty} \frac{f(x)}{x^{\lambda}}\mathrm{d} x = \lim_{a\rightarrow \infty}  \bbint{0}{a} \frac{f(x)}{x^{\lambda}}\mathrm{d} x,
\end{equation}
where it is assumed that $f(x) x^{-\lambda}$ is integrable at infinity so that the limit is guaranteed to exist.

\begin{figure}[t]
	\centering
	\includegraphics[width=.75\textwidth]{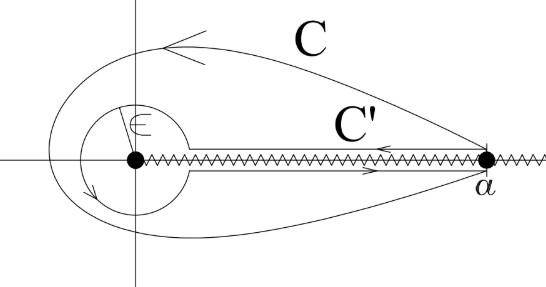}
	\caption{The contour of integration. The contour $C$ does not enclose any poles or branch points of $f(z)$.}
	\label{ch31}
\end{figure}

Now the heart of finite-part integration is the representation of the fundamental finite-part integral as a contour integral in the complex plane \cite{galapon2,galapon1}. The representation requires that $f(x)$ possesses an analytic complex extension, $f(z)$. The contour integral representation depends on the nature of the singularity of the integrand $z^{-\lambda} f(z)$ at the origin, whether $z=0$ is a pole or a branch point.  Assuming that $f(0)\neq 0$ and that $f(z)$ is analytic in the interval $[0,a]$, the contour integral representation of the finite part when $z=0$ is a pole of $z^{-\lambda} f(z)$ is given by 
\begin{equation}\label{polefinitepart}
\bbint{0}{a} \frac{f(x)}{x^{m}} \, \mathrm{d}x = \frac{1}{2 \pi i} \int_{\mathrm{C}} \frac{f(z)}{z^{m}} (\log z - \pi i) \, \mathrm{d}z ,
\end{equation}
for all $\lambda=m = 1,2,3,...$, where $\log z$ is the complex logarithm whose branch cut is the positive real axis and is choosen to coincide with the natural logarithm above the cut, $\mathrm{C}$ is the contour straddling the branch cut of $\log z$ starting from $a$ and ending at $a$ itself, as depicted in Figure-\ref{ch31} \cite{galapon2}. On the other hand, when $z=0$ is a branch point, i.e. $\lambda$ is not an integer, the contour integral representation is given by  
\begin{equation}\label{branchfinitepart}
\bbint{0}{a} \frac{f(x)}{x^{\lambda}} \, \mathrm{d}x = \frac{1}{\mathrm{e}^{-2 \pi \lambda i}-1} \int_{\mathrm{C}} \frac{f(z)}{z^{\lambda} }\, \mathrm{d}z ,
\end{equation} 
where $z^{−\lambda}$ has a branch cut along the positive real axis and is chosen to be positive above the cut, and the contour $\mathrm{C}$ is the same contour as in the representation \eqref{polefinitepart} \cite{galapon2}. In both cases, the contour $\mathrm{C}$ does not enclose any pole of $f(z)$ or intersect any of the branch cuts of $f(z)$. The contour integrals \eqref{polefinitepart} and \eqref{branchfinitepart} are the bridges that connect the two sides of the equality \eqref{byfpi}.

Here we will find that the progenic finite part integrals arise from the singularity of the kernel of transformation.  They are finite parts of divergent integrals of the form 
\begin{equation}
\int_0^d \frac{h(x)}{(d-x)^{\sigma}}\mathrm{d}x,
\end{equation}
 where the divergence arises from a non-integrable singularity at $x=d$. The finite part in this case is obtained by isolating the singularity with the integral 
 \begin{equation}
 \int_0^{d-\epsilon}\frac{h(x)}{(d-x)^{\sigma}}\mathrm{d}x=\tilde{C}_{\epsilon}+\tilde{D}_{\epsilon},
 \end{equation}
 where $\tilde{C}_{\epsilon}$ and $\tilde{D}_{\epsilon}$ are the converging and diverging parts of the integral $\int_c^{d-\epsilon}h(x) (x-d)^{-\sigma} \mathrm{d}x$ as $\epsilon\rightarrow 0$. The finite part integral is then defined as
 \begin{equation}
 \bbint{0}{d} \frac{h(x)}{(d-x)^{\sigma}}\,\mathrm{d}x = \lim_{\epsilon\rightarrow 0} \tilde{C}_{\epsilon}
 \end{equation}
 or equivalently 
 \begin{equation}
 \bbint{0}{d} \frac{h(x)}{(d-x)^{\sigma}}\,\mathrm{d}x = \lim_{\epsilon\rightarrow 0} \left(\int_c^{d-\epsilon} \frac{h(x)}{(d-x)^{\sigma}}\,\mathrm{d}x- \tilde{D}_{\epsilon}\right) .
 \end{equation}
 Now since for $\epsilon>0$ all integrals involved are convergent, we can cast the result in terms of the integral \eqref{divinteg} with the change in variable $x'=d-x$, $\int_c^{d-\epsilon}h(x)\mathrm{d}x=\int_{\epsilon}^{d-c}h(d-x')\mathrm{d}x'$. In this form, everything that we have discussed above for divergent integrals arising from a singularity at the origin holds. 

Two remarks are in order. First, an ambiguity may arise in the definition of the finite part integral. For example, say, after an integration there appeared the term $e^{1/\epsilon}$. This term diverges as $\epsilon\rightarrow 0$. If we take $e^{1/\epsilon}$ as a whole, then the finite part is zero because $\lim_{\epsilon\rightarrow 0}e^{1/\epsilon}=\infty$. However, we can expand it for us to have $e^{1/\epsilon}=1+1/\epsilon + 1/2!\epsilon^2 +\dots$ and now we have the non-zero finite part value equals $1$. Also let us say we arrived at the term $\ln(2 \epsilon)$ which diverges as $\epsilon\rightarrow 0$. Again it is zero when $\ln(2\epsilon)$ is taken as a whole because $\lim_{\epsilon\rightarrow 0} \ln(2\epsilon)=-\infty$; but also $\ln(2\epsilon)=\ln(2) + \ln(\epsilon)$ so that we have the non-zero finite part $\ln(2)$. So which finite part? Here a specific finite part is meant. The finite part is fixed by requiring that the divergent part $D_{\epsilon}$ ($\tilde{D}_{\epsilon}$) only assumes the form of at most the sum of inverse powers of $\epsilon$ and powers of the logarithm $\ln\epsilon$. Then our desired finite parts for the two examples are $1$ and $\ln(2)$, respectively. Our intention for the finite part then requires expansion of $f(x)$ about $x=0$ ($x=d$) in the calculation of the finite part. Here we require that $f(x)$ is infinitely differentiable at the origin so that it admits the representation $f(x)=\sum_{k=0}^{\infty} a_k x^k$ ($\sum_{k=0}^{\infty}b_k (x-b)^k$) in some neighborhood of the origin. 

%%%%%%%%%%%%%%%%%%%%%%%%%%%%%%%%%%%%%%%%%
%%%%%%%%%%%%%%%%%%%%%%%%%%%%%%%%%%%%%%%

Second, we quote from \cite{tica1}: ``from the definition of the finite part integral given by equation \eqref{finitepart}, it is evident that, when  the divergent part $D_{\epsilon}$ vanishes, the finite part integral is just the (Riemann) improper integral. Also the complex contour integral representations of the finite parts given by equations \eqref{polefinitepart} and \eqref{branchfinitepart} reduce to the regular (Riemann) integrals of the integrands when the improper integrals exist. For this reason, we can always replace the regular integral $\int_0^a$ with the finite part integral $\bbint{0}{a}$ without possible confusion for the two values coincide when the former exists. In short, the finite part integral of a convergent integral is just the value of the convergent integral itself.''

Finally, we outline how the method of finite-part integration that will be implemented in this work. Our problem is to evaluate the Stieltjes integral of the form
\begin{equation}\label{convergeint}
\int_0^{\infty} \frac{x^{\lambda-1} f(x)}{(\omega+x)^{\rho}}\mathrm{d}x,\;\;\; \mathrm{Re}(\lambda),\mathrm{Re}(\rho)>0,\, |\mathrm{arg}(\omega)|<\pi,
\end{equation}
where $f(x)$ possesses a complex extension $f(z)$ that is analytic in the integration interval $[0,\infty)$, and with a finite order of zero at the origin. The extension $f(z)$ may have poles and brach points elsewhere in the complex plane. 
 
For the first step of finite-part integration, we induce divergence in the Stieltjes integral \eqref{convergeint} by expanding the kernel of transformation $(\omega+x)^{-\rho}$ about $\omega=0$ as given by \eqref{expand} under the integral, followed by term by term integration. The fundamental divergent integrals are given by
  \begin{equation}\label{divints}
  \int_0^{\infty} \frac{ f(x)}{x^{k+\rho-\lambda+1}}\mathrm{d}x
  \end{equation}
for sufficiently large integer values of $k$. It is assumed in \eqref{divints} that the function under the integral sign is integrable at infinity so that the origin is the only source of non-integrability. This is automatically satisfied under the condition that the Stieltjes integral \eqref{convergeint} exists. The relevant contour integral representation of \eqref{divints} is already given by equation  \eqref{polefinitepart} for the case of pole singularity at the origin, i.e. $\rho-\lambda\in\mathbb{Z}$; and equation \eqref{branchfinitepart} for the case of branch point singularity, i.e. $\rho-\lambda\notin\mathrm{Z}$. The finite part integrals are to be evaluated once the explicit form of $f(x)$ has been specified.
  
 For the second step, we express the Stieltjes integral as a contour integral consistent with the integral representation of the finite part integrals involved. Representing the given integral as complex contour integral requires obtaining the complex extension of the integrand $x^{\lambda-1} f(x)(\omega+x)^{-\rho}$. When it happens that $f(z)$ has at most poles in the complex plane, and $\lambda$ and $\rho$ are integers, the desired complex extension is obtained simply by replacing the real variable $x$ with the complex variable $z$.  However, the extension, which we write as $z^{\lambda-1} f(z)(\omega+z)^{-\rho}$, is generally multivalued, and we choose the branch as follows. When $\lambda$ is a non-integer, the origin is a branch point of $z^{\lambda-1}$ and we choose the branch cut to be $[0,\infty)$, with $z^{\lambda-1}$ positive above the cut. On the other hand, when $\rho$ is a non-integer, the point $z=-\omega$ is a branch point of $(\omega+z)^{-\rho}$ and we choose the branch cut to be $[-\omega,\infty)$, with $(\omega+z)^{-\rho}$ taken to be positive above the cut. Also the extension of $f(x)$, $f(z)$, is choosen such that it is analytic in the integration interval $[0,\infty)$; in addition to that, a crucial requirement on $f(z)$ will be shortly established. These choices define the desired branch for the complex extesion $z^{\lambda-1} f(z)(\omega+z)^{-\rho}$. 
 
 %%%%%%%%%%%%%%%%%
 \begin{figure}[t]
 	\centering
 	\includegraphics[width=0.75\textwidth]{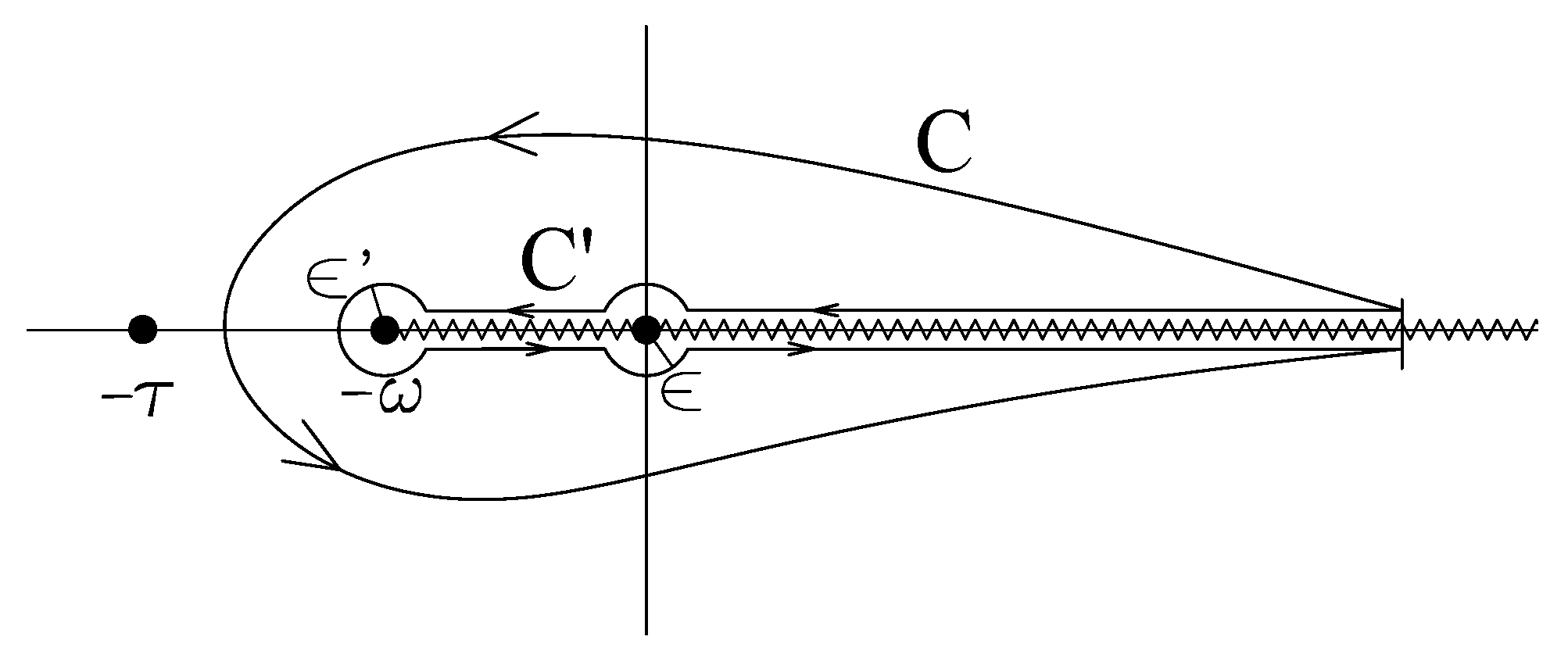}
 	 	\includegraphics[width=0.75\textwidth]{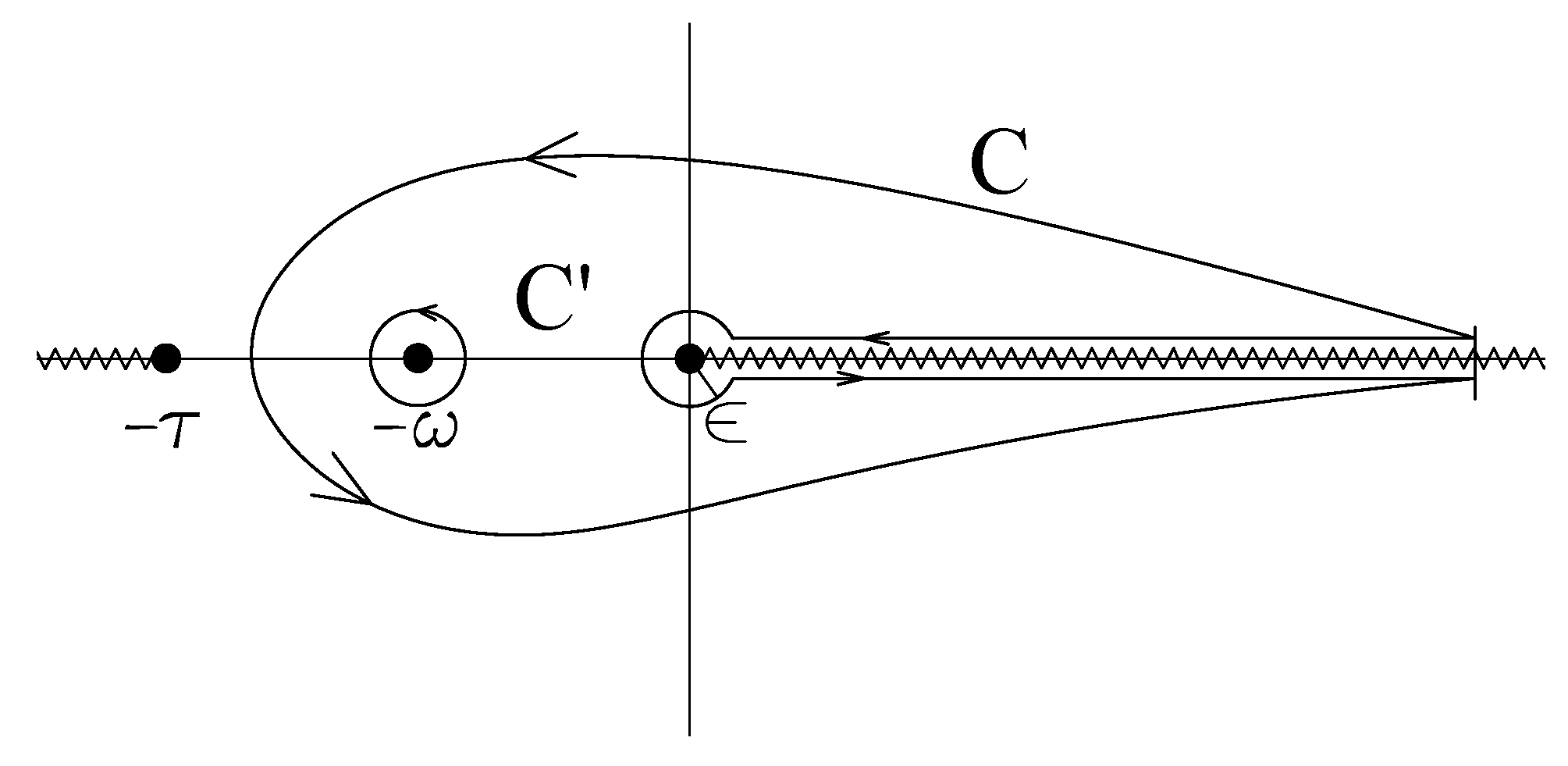}
 	\caption{The deformation of the contour $\mathrm{C}$ to the contour $\mathrm{C}'$ to extract the Stieltjes integral. (Top) The point $z=-\omega$ is a branch point of the kernel. (Bottom) The point $z=-\omega$ is a pole of the kernel of transformation. The point $z=-\tau$ is a branch pole or a pole of the complex extension $f(z)$ of the function $f(x)$.}
 	\label{deformation}
 \end{figure}
 
 Once the extension is fixed,  the Stieltjes integral \eqref{convergeint} is extracted from the contour integral
 \begin{equation}\label{contourint}
 \int_{\mathrm{C}} \frac{z^{\lambda-1} f(z)}{(\omega+z)^{\rho}}\,G(z)\,\mathrm{d}z .
 \end{equation} 
 The contour $\mathrm{C}$ is homotopic to the contour of integration in the representation of the finite part integral. All singularities of $f(z)$, poles or branch points, stay to the right of $\mathrm{C}$ when it is traversed in the positive sense or none of the singularities is enclosed by $\mathrm{C}$ or none of the branch cuts of $f(z)$ intersect $C$. The contour $\mathrm{C}$ is chosen such that the binomial expansion 
 \begin{equation}\label{binomial}
 \frac{1}{(\omega+z)^{\rho}} = \sum_{k=0}^{\infty} {-\rho\choose k} \frac{\omega^k}{z^{k+\rho}}
 \end{equation}
 converges uniformly along the contour $\mathrm{C}$. This requires that the point $z=-b$ must be enclosed by $\mathrm{C}$. The function $G(z)$ is determined by whether $z=0$ is a pole or a branch point of $z^{-(k+\rho-\lambda+1)}$. When the origin is a pole, i.e. when $(\rho-\lambda)\in\mathbb{Z}_0^+$, we have 
 \begin{equation}
 G(z)=\frac{1}{2\pi i}\left(\log z - \pi i\right);  %\;\; z=0\;\; \mathrm{is\; a\; pole\; singularity}
 \end{equation}
 on the other hand, when it is a branch point, i.e. $(\rho-\lambda)\notin\mathbb{Z}$, we have 
 \begin{equation}
     G(z)=\frac{1}{e^{-2\pi (\rho-\lambda) i}-1} . 
 \end{equation}
 
 The desired representation of the Stieltjes integral is obtained from the contour integral \eqref{contourint} by collapsing the contour $\mathrm{C}$ into the real line via the contour $C'$ as depicted in Figure-2. This leads to the representation
 \begin{equation}\label{rep0}
\int_0^{\infty} \frac{x^{\lambda-1} f(x)}{(\omega+x)^{\rho}}\mathrm{d}x = \int_{\mathrm{C}} \frac{z^{\lambda-1} f(z)}{ (\omega+z)^{\rho}}\,G(z)\,\mathrm{d}z + \Delta_{S} , 
 \end{equation}
where $\Delta_{S}$ is a contribution arising from the pole or branch point singularity of the kernel $(\omega+z)^{-\rho}$ at $z=-\omega$,  for which reason we have refered to $\Delta_{S}$ as the singular contribution \cite{tica1,tica2}. In general, when $z=-\omega$ is a branch point, the singular contribution is a progenic finite-part integral. On the other hand, when $z=-\omega$ is a pole, the singular contribution is a residue at $z=-\omega$. Here we consider only the case in which $z=-\omega$ is branch point, except for one illustrative case. 

Now under the assumption that the expansion \eqref{binomial} converges uniformly along the contour $\mathrm{C}$, we can introduce the expansion \eqref{binomial} back into the integral and perform a term by term integration to yield
\begin{equation}\label{bo}
\int_{\mathrm{C}} \frac{f(z)}{z^{\lambda} (\omega+z)^{\rho}}\,G(z)\,\mathrm{d}z = \sum_{k=0}^{\infty} {-\rho \choose k} \omega^k \int_{\mathrm{C}} \frac{f(z)}{z^{k+\rho-\lambda+1}} G(z)\,\mathrm{d}z .
\end{equation}
We recognize that the contour integrals are just the finite parts of the divergent integrals \eqref{divints}. Substituting \eqref{bo} back into equation \eqref{rep0}, we obtain the desired evaluation of the Stieltejes integral in terms of finite part integrals,
\begin{equation}\label{rep1}
\int_0^{\infty} \frac{x^{\lambda-1} f(x)}{(\omega+x)^{\rho}}\mathrm{d}x = \sum_{k=0}^{\infty} {-\rho \choose k} \omega^k \bbint{0}{\infty} \frac{f(x)}{x^{k+\rho-\lambda+1}}\,\mathrm{d}x + \Delta_{S} , 
\end{equation}
where the radius of convergence of the sum is determined by the nearest singularity of $f(z)$ in the complex plane.

The effects of the presence of singularity of $f(z)$ may have gone unnoticed in the above discussion leading to the desired evaluation of the Stieltjes integral in terms of the finite part integral, except, perhaps, the statement following the result \eqref{rep1}. We now highlight them here. The conditions that the contour of representation $\mathrm{C}$ must not intersect any branch cut of $f(z)$ and all singularities of $f(z)$ stay to the left of $\mathrm{C}$ when traversed in the positive sense restrict the possible choices of the complex extension of the function $f(x)$.  These conditions are imposed by our desire to identify the contour integrals in the right hand side of \eqref{bo} as finite-part integrals: the equality in the contour integral representations of the finite part integral given by equations \eqref{polefinitepart} and \eqref{branchfinitepart} hold under the condition that no singularity of $f(z)$ is enclosed by $\mathrm{C}$. These conditions translates to the necessary choice that the branch cut of $f(z)$ cannot overlap with the branch cut of the extension of the kernel of transformation, $(\omega+z)^{-\rho}$.

On the other hand, the condition that the contour $\mathrm{C}$ must enclose the singularity of the kernel, $-\omega$, together with the condition that the singularities of $f(z)$ must be to the left of $\mathrm{C}$, requires that  $-\omega$ cannot be farther from the origin than the nearest singularity of $f(z)$. This implies that the radius of convergence of the infinite series in the right hand side of \eqref{rep1} cannot be larger than the distance of the nearest singularity of $f(z)$ to the origin. In contrast, under the condition that $f(z)$ is entire \cite{galapon2,tica1,tica2}, $f(z)$ has no singularity in the finite complex plane so that the radius of convergence of the corresponding series of finite part integrals is infinite. An unexpected consequence of our considerations is that, once all conditions have been satisfied, the nature of the singularities of $f(z)$, be it a pole or a branch point, does not have a bearing in what follows to arrive at the finite part integration \eqref{rep1} of the Stieltjes integral. The singularities are manifested only by imposing restriction of the radius of convergence of the infinite series of finite part integrals. 

\section{The Fundamental Finite-Part Integrals} \label{fundamental}
We now identify and evaluate the fundamental finite part integrals arising from finite-part integration of the integral representations \eqref{mainlemma} and \eqref{general}. For the reprsentation \eqref{mainlemma}, we identify the function $f(x)$ in \eqref{convergeint} to be $(a+x)^{-\mu}$ with the corresponding kernel $(b+x)^{-\rho}$. Then, from the general expression for the fundamental divergent integrals for the Stieltjes transform given by equation \eqref{divints}, the fundamental divergent integrals for the evaluation of \eqref{mainlemma} are specific cases of the divergent integral  
\begin{equation}\label{funint}
\int_0^{\infty} \frac{(s+t)^{-\upsilon}}{t^{\lambda}}\mathrm{d}t, 
\end{equation}
with $\mathrm{Re}(\lambda)\geq 1$, $\upsilon\neq 0$, $\mathrm{Re}(\lambda+\upsilon)>1$ and $|\mathrm{arg}(s)|<\pi$, where the second and third conditions ensure that the origin is the only source of divergence of the integral. We shall find in Section-\ref{section3f2} that the fundamental divergent integrals for the case of the generalized hypergeometric function $_3F_2$ for the specific set of parameters considered there are likewise specific cases of \eqref{funint}. The divergence of \eqref{funint} is due to a branch point singularity at the origin when $\lambda$ is not an integer, and is due to a pole singularity when $\lambda$ is a positive integer. As the contour integral representations of the finite part integral already make apparent, the finite part depends on the type of non-integrable singularity at the origin of the divergent integral.

\subsection{Branch point singularity}
We now obtain the finite part of the divergent integral \eqref{funint} for non-integer $\lambda$, specifically $\mathrm{Re}(\lambda)\notin\mathbb{Z}_0^+$.	Since the upper limit of integration extends to infinity, the desired finite part is given by equation \eqref{finitepart3}. We then cut the upper limit to some finite $c>0$, compute the finite part of the divergent integral $\int_0^c \mathrm{d}x\, (s+x)^{-\mu} x^{-\lambda}$, and eventually take the limit $c\rightarrow\infty$. To facilitate the computation of the finite-part, we temporarily impose the condition $c<|s|$, with the intention to eventually use analytic continuation to go beyond $c>|s|$. Then  for some arbitrarily small positive $\epsilon$, we have
	\begin{equation}
	\int_{\epsilon}^{c} \frac{(s+x)^{-\upsilon}}{x^{\lambda}} \,\mathrm{d}x = \sum_{l=0}^{\infty} {-\upsilon \choose l}  \frac{1}{s^{l+\upsilon}} \, \int_{\epsilon}^{c} x^{l-\lambda} \,\mathrm{d}x,
	\label{FPIcase2}
	\end{equation}
	where the right hand side is obtained by expanding $(a+x)^{-\upsilon}$ about $x=0$ and distributing the integration. The term by term integration in equation \eqref{FPIcase2} is allowed under the assumed condition $c<|s|$. 
	
	Let $\lambda=m+\sigma$, where $m=1, 2, \dots$ is the integer part of the real part of $\lambda$ and $\sigma\notin\mathbb{Z}$, in particular $\mathrm{Re}(\sigma)>0$. Then the summation in \eqref{FPIcase2} can be divided into two groups of terms in $\epsilon$: the terms with negative exponents corresponding to the indices $0 \leq l \leq m -1$;  and the terms with positive exponents corresponding to the indices $ m \leq l \leq \infty$. We have
	\begin{equation}
	\label{FPIc2.1}
	\begin{split}
	\int_{\epsilon}^{c} \frac{(s+x)^{-\upsilon}}{x^{\lambda}}  \, \mathrm{d}x
	=&\sum_{l=0}^{m-1} \frac{1}{s^{l+\upsilon}} {-\upsilon \choose l} \,   \,\left( \frac{1}{c^{m+\sigma+1-l}} - \frac{1}{\epsilon^{m+\sigma+1-l}} \right) \\
	& + \sum_{l=m}^{\infty} \frac{1}{s^{l+\upsilon}} {-\upsilon \choose l}\, \frac{1}{(l-m-\sigma+1)}  \, \left( c^{l-m-\sigma+1}-\epsilon^{l-m-\sigma+1} \right).
	\end{split}
	\end{equation}
	The converging and diverging parts of the left hand side of equation \eqref{FPIc2.1} can now be extracted from the right hand side of the equality. They are respectively given by
	\begin{equation}
	\begin{split}
	C_{\epsilon} =& \sum_{l=0}^{m-1} \frac{1}{s^{l+\upsilon}} {-\upsilon \choose l} \, \frac{1}{(l-m-\sigma+1)}  \, \frac{1}{c^{m+\sigma+1-l}} \\
	& + \sum_{l=m}^{\infty} \frac{1}{s^{l+\upsilon}} {-\upsilon \choose l} \, \frac{1}{(l-m-\sigma+1)}  \, \left( c^{l-m-\sigma+1}-\epsilon^{l-m-\sigma+1} \right)
	\end{split}
	\end{equation}
	\begin{equation}
	D_{\epsilon} = \sum_{l=0}^{m-1} \frac{1}{s^{l+\upsilon}} {-\upsilon \choose l} \, \frac{1}{(l-m-\sigma+1)}  \, \frac{1}{\epsilon^{m+\sigma+1-l}}
	\end{equation}
	Then the finite-part is obtained by dropping the contribution of the diverging part and assigning the value of the converging part in the limit $\epsilon\rightarrow 0$ as the value of the divergent integral,
	\begin{equation}
	\begin{split}
	\bbint{0}{c} \frac{(s+x)^{-\upsilon}}{x^{m+\sigma}} \,\mathrm{d}x = &  \sum_{l=0}^{m-1} \frac{1}{s^{l+\upsilon}} {-\upsilon \choose l} \, \frac{1}{(l-m-\sigma+1)}  \, \frac{1}{c^{m+\sigma+1-l}} \\
	& + \sum_{l=m}^{\infty} \frac{1}{s^{l+\mu}} {-\upsilon \choose l} \, \frac{1}{(l-m-\sigma+1)}  \, c^{l-m-\sigma+1} .
	\end{split}
	\end{equation}
	
	We now take the limit as $c \to \infty$. The first term  extends to all $c$ and vanishes in the limit. However, the second term is valid only for $c<|s|$. We sum it and extend the sum analytically everywhere. The sum is given by 
	\begin{equation}\label{FPIc2v2.1}
	\begin{split}
	\sum_{l=m}^{\infty} \frac{1}{s^{l+\upsilon}} {-\upsilon \choose l} \, \frac{c^{l-m-\sigma+1}}{(l-m-\sigma+1)} = \frac{c^{1-\sigma}}{s^{\upsilon+m}} \frac{1}{(1-\sigma)} {-\upsilon \choose m}
	\pFq{3}{2}{1, m +\upsilon, 1-\sigma}{m+1, 2-\sigma}{-\frac{c}{s}} .
	\end{split}
	\end{equation}
	The right hand side is now valid for all values of $a$ and $c$. Then the finite part is given by the limit
	\begin{equation}\label{limit}
	\bbint{0}{\infty} \frac{(s+x)^{-\upsilon}}{x^{m+\sigma}}\,\mathrm{d}x = \lim_{c\rightarrow\infty} \frac{c^{1-\sigma}}{s^{\upsilon+m}} \frac{1}{(1-\sigma)} {-\upsilon \choose m}\\
	\pFq{3}{2}{1, m +\upsilon, 1-\sigma}{m+1, 2-\sigma}{-\frac{c}{s}} .
	\end{equation}
	
	The limit is obtained by means of the asymptotic expansion of the hypergeometric function $_3F_2$ for large arguments \cite{simplepole3F2},
	\begin{equation}
	\label{3F2simple} 
	\begin{split}
	&\pFq{3}{2}{a_1, a_2, a_3}{b_1, b_2}{z} =\frac{\Gamma(b_1) \Gamma(b_2)}{\Gamma(a_1) \Gamma(a_2) \Gamma(a_3)} \\
	&\hspace{16mm}\times \left[ \frac{\Gamma(a_1) \Gamma(a_2-a_1) \Gamma(a_3-a_1)}{\Gamma(b_1-a_1) \Gamma(b_2-a_1) \Gamma(a_3)} (-z)^{-a_1} \left(  1+\mathcal{O} \left( \frac{1}{z} \right) \right) \right.\\ 
	&\hspace{16mm}\left. 
	+ \frac{\Gamma(a_2) \Gamma(a_1-a_2) \Gamma(a_3-a_2)}{\Gamma(b_1-a_2) \Gamma(b_2-a_2) \Gamma(a_3)} (-z)^{-a_2} \left(  1+\mathcal{O} \left( \frac{1}{z} \right) \right) \right.\\ 
	&\hspace{16mm}\left. +  \frac{\Gamma(a_3) \Gamma(a_1-a_3) \Gamma(a_2-a_3)}{\Gamma(b_1-a_3) \Gamma(b_2-a_3) \Gamma(a_3)} (-z)^{-a_3} \left(  1+\mathcal{O} \left( \frac{1}{z} \right) \right) \right],  |z| \to \infty.
	\end{split}
	\end{equation}
	Applying this asymptotic expansion to the hypergeometric function in \eqref{limit} and performing some simplifications, the finite-part integral \eqref{limit} assumes the form
	\begin{equation}\label{rawfinitepart}
	\begin{split}
	&\bbint{0}{\infty} \frac{(s+x)^{-\upsilon}}{x^{m+\sigma}} \, \mathrm{d}x \\
	&\hspace{8mm}= \lim_{c \to \infty} \, \frac{c^{1-\sigma}}{s^{\upsilon+m}} {-\upsilon \choose m} \frac{\Gamma(m+1)}{\Gamma(m+\upsilon)}\Bigg[ \frac{\Gamma(m+\upsilon-1) \Gamma(-\sigma)}{\Gamma(m) \Gamma(1-\sigma) } \left(\frac{s}{c}\right) \left(1+O\left(\frac{1}{c}\right)\right)\\
	&\hspace{20mm} + \frac{\Gamma(m+\upsilon) \Gamma(1-m-\upsilon) \Gamma(1-\sigma-m-\upsilon)}{\Gamma(1-\upsilon) \Gamma(2-\sigma-m-\upsilon)} \left(\frac{s}{c}\right)^{m+\upsilon}\left(1+O\left(\frac{1}{c}\right)\right)\\
	& \hspace{20mm} + \frac{\Gamma(1-\sigma) \Gamma(\sigma) \Gamma(m+\upsilon-1+\sigma)}{\Gamma(\sigma)} \left(\frac{s}{c}\right)^{1-\sigma} \left(1+O\left(\frac{1}{c}\right)\right)\Bigg] .
	\end{split}
	\end{equation}
	
	Since $\mathrm{Re}(\sigma)>0$, the first term in the right hand side of \eqref{rawfinitepart} vanishes in the limit as $c\rightarrow\infty$; likewise, the second term vanishes in the same limit. Only the third term contributes in the limit. Then the finite part becomes
	\begin{equation}\label{fpix}
	\bbint{0}{\infty} \frac{(s+x)^{-\upsilon}}{x^{m+\sigma}} \, \mathrm{d}x	= \frac{(-1)^m}{s^{m+\sigma+\upsilon-1}} \frac{\Gamma(\sigma)\Gamma(1-\sigma) \Gamma(m+\sigma+\upsilon-1)}{ \Gamma(m+\sigma)\Gamma(\upsilon)} ,
	\end{equation}
	where we have used the following identities to perform the simplifications leading to this result,
	\begin{equation}\label{identity}
	\Gamma(a+k)=(a)_k \, \Gamma(a), \;\;\; {n \choose k} = (-1)^k \frac{(-n)_k}{k!}
	\end{equation}
	We wish to express equation \eqref{fpix} in terms of the original parameter $\lambda$. This is accomplished using the identities $\Gamma(\sigma)\Gamma(1-\sigma)=\pi \csc(\pi \sigma)$ and $\sin(\pi\sigma)= (-1)^m \sin(\pi(m+\sigma))$. The finite part becomes
	\begin{equation}\label{another}
	\bbint{0}{\infty} \frac{(s+x)^{-\upsilon}}{x^{m+\sigma}} \, \mathrm{d}x	= \frac{\pi\, \Gamma(m+\sigma+\upsilon-1)}{s^{m+\sigma+\upsilon-1}\, \sin(\pi(m+\sigma)) \Gamma(m+\sigma)\Gamma(\upsilon)}.
	\end{equation}
	With $\lambda=m+\sigma$, this reduces to the desired form of the finite part integral,
	\begin{equation}\label{finxx}
	\bbint{0}{\infty} \frac{(s+x)^{-\upsilon}}{x^{\lambda}} \, \mathrm{d}x = \frac{\pi\, \Gamma(\lambda+\upsilon-1)}{s^{\lambda+\upsilon-1}\, \sin(\pi\lambda)\, \Gamma(\upsilon)\Gamma(\lambda)},
	\end{equation} 
	for $\mathrm{Re}(\lambda)>1$, $\nu\neq 0$, $\mathrm{Re}(\lambda)\notin\mathbb{Z}^+$, $\mathrm{Re}(\lambda+\upsilon)>1$, and $|\mathrm{arg}(s)|<\pi$. 

When $0<\mathrm{Re}(\lambda)<1$ (and for sufficiently large $\mathrm{Re}(\upsilon)$), the integral \eqref{funint} converges. To evaluate the (this time convergent) integral, we use the known integral \cite[3.194.3,p-318]{gr}
\begin{equation}
\int_0^{\infty} \frac{x^{d-1}}{(1+c x)^{\sigma}}\,\mathrm{d}x = c^{-d} B(d,\sigma-d),\;\; |\mathrm{arg}\,c|<\pi , \;\; \mathrm{Re}\, \sigma > \mathrm{Re}\, d>0,
\end{equation}
where $B(a,b)=\Gamma(a)\Gamma(b)/\Gamma(a+b)$ is the beta function. Applying this result to integral \eqref{funint}, we obtain the result
\begin{equation}\label{finiteres}
\int_0^{\infty} \frac{(s+x)^{-\upsilon}}{ x^{\lambda}} \,\mathrm{d}x
= \frac{\Gamma(1-\lambda)\Gamma(\upsilon-1+\lambda)}{s^{\upsilon+\lambda-1}\Gamma(\upsilon)}
\end{equation}
 for $\mathrm{Re}(\upsilon)>\mathrm{Re}(1-\lambda)>0$. We can recast this expression further using the reflection formula $\Gamma(-z)=-\pi/\sin(\pi z) \Gamma(1+z)$ to rewrite the factor $\Gamma(1-\lambda)$. The result is
\begin{equation}
\int_0^{\infty} \frac{(s+x)^{-\upsilon}}{ x^{\lambda}} \,\mathrm{d}x
= \frac{\pi\, \Gamma(\lambda+\upsilon-1)}{s^{\lambda+\upsilon-1}\, \sin(\pi\lambda)\, \Gamma(\upsilon)\Gamma(\lambda)}, 
\label{BetaInt}
\end{equation}
for $0<\mathrm{Re}(\lambda)<1$ and $\mathrm{Re}(\upsilon+\lambda)>1$. Comparing this with the finite part \eqref{finxx}, we find that the finite part integral is the extension of the right hand side of \eqref{BetaInt} outside its strip of analyticity.  

Conversely, observe that when the integral converges, examination of equation \eqref{FPIcase2} shows that the divergent part vanishes in the limit $\epsilon\rightarrow 0$. This vanishing of the divergent term is equivalent to just dropping the term so that the value of the regular convergent integral coincides with the value of the finite part integral \eqref{finxx} extended to the value of $\lambda$ that makes the integral convergent. This just demonstrates the fact that the finite part of a convergent integral coincides with the regular value of the integral. This allows us to combine \eqref{finxx} and \eqref{BetaInt} into the single expression which is just \eqref{finxx} itself, and extend the domain of the finite part integral \eqref{finxx} to all $\mathrm{Re}(\lambda)>0$ with all the other conditions unchanged.

\subsection{Pole singularity}
We now obtain the finite part when $\lambda$ is an integer, i.e. $\lambda=n+1$ for all $n\in\mathbb{Z}^+_0$.	We proceed in the same way we did above for the evaluation of the finite part. Then 
\begin{equation}
\bbint{0}{\infty} \frac{(s+x)^{-\upsilon}}{x^{n+1}} \, \mathrm{d}x = \lim_{c\rightarrow\infty}\bbint{0}{c} \frac{(s+x)^{-\upsilon}}{x^{n+1}} \, \mathrm{d}x ,
\end{equation}
for $\mathrm{Re}(n+\upsilon)>0$. We now extract for the finite part for some finite $c$. Let $c$ be sufficiently small so that we can expand the integrand and perform a term by term integration, that is, 
\begin{equation}
\int_{\epsilon}^{c} \frac{(s+x)^{-\upsilon}}{x^{n+1}} \,\mathrm{d}x =  \sum_{l=0}^{\infty} {-\upsilon \choose l}  \frac{1}{s^{l+\upsilon}} \, \int_{\epsilon}^{c} x^{l-n-1} \,\mathrm{d}x .
\label{FPIcase1}
\end{equation}
We split the summation into three parts: $0 \leq l \leq n-1 $, $l = n$, and $n+1 \leq l \leq \infty$. Performing the indicated integrations, we obtain
\begin{equation}
\begin{split}
& \int_{\epsilon}^{c} \frac{(s+x)^{-\upsilon}}{x^{n+1}} \, \mathrm{d}x =  \sum_{l=0}^{n-1} {-\upsilon \choose l}\frac{1}{s^{l+\upsilon}} \, \frac{1}{(l-n)} \, \left( \frac{1}{c^{-l+n}} - \frac{1}{\epsilon^{-l+n}} \right) \\& \qquad + {-\upsilon \choose n}\frac{1}{s^{n+\mu}} \left( \ln{c} - \ln{\epsilon} \right) + \sum_{l=n+1}^{\infty}{-\upsilon \choose l} \frac{1}{s^{l+\upsilon}} \, \frac{1}{(l-n)} \, \left( c^{l-n}-\epsilon^{l-n} \right)
\label{FPIc1} .
\end{split}
\end{equation}
From this expression, we extract the converging part in the limit of arbirary small $\epsilon$, 
\begin{equation}
\begin{split}
C_{\epsilon} = & \sum_{l=0}^{k+\alpha-1} {-\upsilon \choose l}\frac{1}{s^{l+\upsilon}} \, \frac{1}{(l-n)}  \, \frac{1}{c^{-l+n}} +{-\upsilon \choose n} \frac{\ln{c}}{s^{l+n+\upsilon}}\\  
&\hspace{30mm}+ \sum_{l=n+1}^{\infty} {-\upsilon \choose l}\frac{1}{s^{l+\upsilon}} \, \frac{1}{(l-n)}  \, \left( c^{l-n}-\epsilon^{l-n} \right) .
\end{split}
\end{equation}
	
Then from the definition of the finite part, we obtain after taking the limit of the converging part, 
\begin{equation}
\begin{split}
\bbint{0}{c} &\frac{(s+x)^{-\upsilon}}{x^{n+1}}  \,\mathrm{d}x =  \sum_{l=0}^{n-1} {-\upsilon \choose l}\frac{1}{s^{l+\upsilon}} \, \frac{1}{(l-n)}  \, \frac{1}{c^{-l+n}} \\
&\hspace{20mm}+ {-\upsilon \choose n}\frac{1}{s^{n+\upsilon}} \, \ln{c}  + \sum_{l=n+1}^{\infty} {-\upsilon \choose l}\frac{1}{s^{l+\upsilon}} \, \frac{1}{(l-n)}  \, c^{l-n}
\end{split}
\end{equation}
Moreover, by taking $c \to \infty$, the finite-part integral reduces to
\begin{equation}
\begin{split}
\bbint{0}{\infty} \frac{(s+x)^{-\upsilon}}{x^{n+1}} \, \mathrm{d}x = \lim_{c \to \infty} \left[  \frac{1}{s^{n+\upsilon}} {-\upsilon \choose n}\, \ln{c} + \sum_{l=n+1}^{\infty} \frac{1}{s^{l+\mu}} {-\upsilon \choose l}\, \frac{1}{(l-n)}  \, c^{l-n} \right].
\label{FPIc1v2}
\end{split}
\end{equation}
To implement the limit, we sum the infinite series by hypergeometric summation. The result is
\begin{equation}
\sum_{l=n+1}^{\infty} \frac{1}{s^{l+\upsilon}}{-\upsilon \choose l} \, \frac{c^{l-n}}{(l-n)}  = {-\upsilon \choose n+1} \frac{c}{s^{n+\upsilon+1}} \, \pFq{3}{2}{1,1,n+\mu+1}{2,n+2}{-\frac{c}{s}}.
\end{equation}

The asymptotic expansion of the hypergeometric function $_3F_2$ for the case of double pole with large argument \cite{doublepole3F2} is given by
\begin{equation}
\begin{split}
&\pFq{3}{2}{a_1,a_1,a_3}{b_1,b_2}{z} 
= \frac{\Gamma(b_1) \Gamma(b_2) \Gamma^2(a_1-a_3)} {\Gamma^2(a_1) \Gamma(b_1-a_3) \Gamma(b_2-a_3)} (-z)^{-a_3} \left( 1+\mathcal{O} \left( \frac{1}{z} \right) \right) \\
&\hspace{6mm} + \frac{\Gamma(b_1) \Gamma(b_2) \Gamma(a_3-a_1)}{\Gamma(a_1) \Gamma(a_3) \Gamma(b_1-a_1) \Gamma(b_2-a_1)} \, \Psi(z) (-z)^{-a_1} \left( 1+\mathcal{O} \left( \frac{1}{z} \right) \right), |z| \to \infty,
\end{split}
\label{3F2doublepole}
\end{equation}
where the function $\Psi(z)$ is given by
\begin{equation}
\Psi(z) = \log(-z)+\psi(a_3-a_1)-\psi(b_1-a_1)-\psi(b_2-a_1)-\psi(a_1)-2\gamma, 
\end{equation}
in which $\psi(x)$ is a digamma function and $\gamma = -\psi(1)$ is the Euler–Mascheroni constant. By the virtue of \eqref{3F2doublepole}, the second term of \eqref{FPIc1v2} becomes
\begin{equation}
\begin{split}
& \frac{c}{s^{n+\upsilon+1}} \, \pFq{3}{2}{1,1,n+\upsilon+1}{2,n+2}{-\frac{c}{s}}  = \frac{\Gamma(n+2) \Gamma^2(-n-\upsilon)} {\Gamma(1-\upsilon) \Gamma(1-n-\upsilon)} \frac{1}{c^{n+\upsilon}} \left(1+\mathcal{O}\left(\frac{s}{c}\right)\right)\\
&\hspace{8mm} + \frac{1}{s^{n+\upsilon}} \frac{\Gamma(n+2) \Gamma(n+\upsilon)}{\Gamma(n+\upsilon+1) \Gamma(n+1)} \left(\ln\left(\frac{c}{s}\right) + \psi(n+\upsilon)-\psi(n+1) \right)
\left(1+\mathcal{O}\left(\frac{s}{c}\right)\right),
\label{3f2dp}
\end{split}
\end{equation}
as $c\rightarrow\infty$. For $\upsilon$ not equal to a positive integer, the first term in the right hand  side of \eqref{3f2dp} vanishes in the limit $c\rightarrow\infty$ since $(n+\upsilon)>0$; provided we interpret the first term as a limit when $\nu$ is a positive integer, the first term vanishes as well when $\nu$ is a positive integer. On the other hand, the second term contributes a term proportional to $\ln c$ which is expected to cancel the first term in equation \eqref{FPIc1v2}. Indeed the logarithmic term $\ln(c)$ cancels out with the fact that
\begin{equation}
{-\upsilon \choose n+1} \frac{\Gamma(n+2)  \Gamma(n+\upsilon)}{\Gamma(n+\upsilon+1)\Gamma(n+1)} = -{-\upsilon \choose n}.
\end{equation} Substituting \eqref{3f2dp} back into \eqref{FPIc1v2} and performing the limit, we obtain the desired finite-part integral,
\begin{equation}
\bbint{0}{\infty} \frac{(s+x)^{-\upsilon}}{x^{n+1}} \, \mathrm{d}x = \frac{(-1)^n (\upsilon)_n}{s^{n+\upsilon} \,n!} \left(\ln(s) + \psi(n+1) - \psi(n+\upsilon) \right),
\label{FPIpole}
\end{equation}
for all $n\in\mathbb{Z}^+_0$, $\upsilon\neq 0$, $\mathrm{Re}(n+\upsilon)>0$, and $|\mathrm{arg}(s)|<\pi$.

\section{Transformation equations arising for the Gauss function $_2F_1$}\label{origin}
We now work out the consequences of the integral representation \eqref{mainlemma} for the Gauss function 
arising from direct finite-part integration of the  Stieltjes integral
\begin{equation}\label{si1}
\int_0^{\infty} \frac{x^{\nu-1} }{(a+x)^{\mu} (b+x)^{\rho}} \, \mathrm{d}x,
\end{equation}
for $\mathrm{Re}(\nu)>0$, $\mathrm{Re}(\mu+\rho-\nu)>0$, $|\mathrm{arg}(a)|<\pi$ and $|\mathrm{arg}(b)|<\pi$. We have chosen to interpret this as a Stieltjes transform of the function $x^{\nu-1} (a+x)^{-\mu}$ with the the kernel of transformation $(b+x)^{-\rho}$. With $(b+x)^{-\rho}$ as the kernel, we require that $\mathrm{Re}(\rho)>0$. The function $f(x)$ in \eqref{convergeint} corresponds to the function $(a+x)^{-\mu}$. We will temporarily impose the restriction that $a>b>0$ and lift this restriction later by appealing to analytic continuation to extend the result in the complex plane. The requirement $b<a$ is necessary to satisfy the condition that the nearest singularity of the complex extension of $f(x)$ must be farther from the origin than the singularity of the kernel. The complex extension of the kernel is $(b+z)^{-\rho}$ and has a pole or branch point singularity at $-b$ depending on whether $\rho$ is a positive integer or not. When $\rho$ is a positive non-integer, $-b$ is a branch point and we choose its branch cut to be $[-b,\infty)$, in accordance with the general requirements stated in Section-\ref{finitepartintegration}. The complex extension, $(a+z)^{-\mu}$, of $(a+x)^{-\mu}$ has $-a$ as its lone singularity. If $\mu$ is a positive integer, the singularity is a pole; otherwise, it is a branch point. If $-a$ happens to be a branch point, the branch cut is chosen such that the contour $\mathrm{C}$ cannot cross the cut. We choose the branch cut $(-\infty,-a]$. 

The fundamental finite part integrals for the finite part integration of \eqref{si1} are explicitly given by
\begin{equation}
\int_0^{\infty} \frac{(a+x)^{-\mu}}{x^{k+\rho-\nu+1}}\,\mathrm{d}x,
\end{equation}
for sufficiently large positive integer values of $k$, where $\mathrm{Re}(\rho+\mu-\nu)>0$ to guarantee that the divergence arises from the origin only. Since $(k+1)$ is already an integer, there are two possible cases depending on whether the difference $(\rho-\nu)$ is an integer or not, and each case further splits according to the sign of the difference.  

\subsection{Case $(\rho-\nu)\notin\mathbb{Z}$}\label{case3} We will consider separately the cases $\rho\notin\mathbb{Z}$ and $\rho\in\mathbb{Z}^+$ corresponding to $-b$ being a branch point and a pole of the kernel, respectively. For a reason to be appreciated below, it is only for the present case that we will consider the case when the kernel has a pole singularity.
  
\subsubsection{Case $\rho\notin\mathbb{Z}$}For $(\rho-\nu)\notin\mathbb{Z}$, the divergence arises from a branch point singularity at the origin. We extract the Stieltjes integral from the contour integral 
\begin{equation}
 \frac{1}{e^{-2\pi i(\rho-\nu)}-1} \int_{\mathrm{C}} \frac{z^{\nu-1} (a+z)^{-\mu}}{(b+z)^{\rho}} \, \mathrm{d}z,
\end{equation}
where C is the contour in Figure-\ref{deformation} where $\omega=b$ and $\tau=a$. Under the condition that $a>b$, whether $z=-a$ is a pole or branch point of $(a+z)^{-\mu}$ does not matter. To extract the Stieltjes integral from the above integral, we deform the contour $\mathrm{C}$ into the contour $\mathrm{C}'$ as shown in Figure-\ref{deformation} (top figure). We obtain the following representation of the given integral, 
\begin{equation}\label{meme}
\begin{split}
&\int_{0}^{\infty}  \frac{x^{\nu-1} \, (a+x)^{-\mu}}{(b+x)^{\rho}} \, \mathrm{d}x
= \frac{1}{e^{-2\pi i (\rho-\nu)}-1} \int_C \frac{z^{\nu-1}(a+z)^{-\mu}}{(b+z)^{\rho}} \,\mathrm{d}z \\
&+\frac{\sin(\pi\rho)}{\sin(\pi(\rho-\nu))} \lim_{\epsilon\rightarrow 0}\left[\int_0^{b-\epsilon} \frac{x^{\nu-1}(a-x)^{-\mu}}{(b-x)^{\rho}}\mathrm{d}x+ \frac{e^{i\pi(\rho-\nu)}}{2i \sin(\pi\rho)} \int_{C_{\epsilon}}\frac{z^{\nu-1}(a+z)^{-\mu}}{(b+z)^{\rho}} \,\mathrm{d}z\right].
 \end{split}
\end{equation}
In a short while we will show that the limit is the finite part of the divergent integral $\int_0^b x^{\nu-1} (a-x)^{-\mu} (b-x)^{-\rho}\,\mathrm{d}x$. Now under the condition imposed upon the contour $\mathrm{C}$, we can expand the kernel $(b+z)^{-\rho}$ in powers of $b/z$ and perform a term by term integration. We obtain
\begin{equation}
\begin{split}\label{sese}
\int_{0}^{\infty}  \frac{x^{\nu-1} \, (a+x)^{-\mu}}{(b+x)^{\rho}} \, \mathrm{d}x & =   \sum_{k=0}^{\infty} {- \rho \choose k} \, b^{k} \, \bbint{0}{\infty} \frac{(a+x)^{-\mu}}{x^{k+\rho-\nu+1}} \,\mathrm{d}x
\\
&\hspace{12mm}+\frac{\sin(\pi\rho)}{\sin(\pi(\rho-\nu))} \bbint{0}{b} \frac{x^{\nu-1}(a-x)^{-\mu}}{(b-x)^{\rho}}\mathrm{d}x,
\end{split}
\end{equation}
which is valid for $a>b>0$, $\mu\neq 0$, $\mathrm{Re}(\nu)>0$, $(\rho-\nu)\notin\mathbb{Z}$, $\rho\notin\mathbb{Z}$ and $\mathrm{Re}(\rho+\mu-\nu)>0$. Equation \eqref{sese} is the desired evaluation in terms of the finite part integral and the singular contribution, which is a progenic finite-part integral arising from the singularity of the kernel at $z=-b$. 

We now show that the limit in \eqref{meme} is indeed a finite part integral. To analyze the bracketed quantity in the second term, we evaluate the integral around the circular contour $\mathrm{C}_{\epsilon}$. We parameterize the variable of integration by $z=-b + \epsilon e^{i\theta}$, with $0<\theta\leq 2\pi$, and perform a binomial expansion of $(-b+\epsilon e^{i\theta})^{\nu-1}$ and $(a-b+\epsilon e^{i\theta })^{-\mu}$. For sufficiently small $\epsilon$,
the integration can be distributed in the resulting double summation. Upon collecting equal powers of $\epsilon$ and performing the indicated integration, we obtain the result
\begin{equation}\label{bebe}
\begin{split}
&\frac{e^{i\pi(\rho-\nu)}}{2 i \sin(\pi \rho)} \int_{\mathrm{C}_{\epsilon}}
\frac{z^{\nu-1} (a+z)^{-\mu}}{(b+z)^{\rho}}\,\mathrm{d}z\\
&\hspace{8mm} = \frac{b^{\nu-1}}{(a-b)^{\mu}} \sum_{k=0}^{\infty}
\frac{(-1)^k}{b^k} \left[\sum_{l=0}^k (-1)^l {\nu-1\choose k-l} {-\mu \choose l}
\frac{b^l}{(a-b)^l)}\right] \frac{\epsilon^{k-\rho+1}}{(k-\rho+1)}.
\end{split}
\end{equation} 
Observe that the first $(\lfloor\rho\rfloor-1)$ terms diverge as $\epsilon\rightarrow 0$. These terms exactly cancel the diverging terms of the integral $\int_0^{b-\epsilon} x^{\nu-1} (a-x)^{-\mu} (b-x)^{-\rho}\,\mathrm{d}x$ in the same limit. To see this let us consider a positive number $c$ less than but sufficiently close to $(b-\epsilon)$, and split the integral into two,   
\begin{equation}
\begin{split}
\int_0^{b-\epsilon} \frac{x^{\nu-1}(a-x)^{-\mu}}{(b-x)^{\rho}}\mathrm{d}x = \int_0^{c} \frac{x^{\nu-1}(a-x)^{-\mu}}{(b-x)^{\rho}}\mathrm{d}x + \int_c^{b-\epsilon} \frac{x^{\nu-1}(a-x)^{-\mu}}{(b-x)^{\rho}}\mathrm{d}x .
\end{split}
\end{equation} 
Since $x^{\nu-1} (a-x)^{-\mu}$ is analytic at $x=b$, we perform a binomial expansion of it at $x=b$. We can always choose $c$ such that it falls within the region of convergence of the expansion. Collecting equal powers of $x$ and distributing the integration, we obtain 
\begin{equation}\label{boy}
\begin{split}
&\int_0^{b-\epsilon} \frac{x^{\nu-1}(a-x)^{-\mu}}{(b-x)^{\rho}}\mathrm{d}x = \int_0^{c} \frac{x^{\nu-1}(a-x)^{-\mu}}{(b-x)^{\rho}}\mathrm{d}x \\
&\hspace{12mm}+\frac{b^{\nu-1}}{(a-b)^{\mu}} \sum_{k=0}^{\infty}
\frac{(-1)^k}{b^k} \left[\sum_{l=0}^k (-1)^l {\nu-1\choose k-l} {-\mu \choose l}
\frac{b^l}{(a-b)^l)}\right] \frac{(b-c)^{k-\rho+1}}{(k-\rho+1)}\\
&\hspace{12mm}-\frac{b^{\nu-1}}{(a-b)^{\mu}} \sum_{k=0}^{\infty}
\frac{(-1)^k}{b^k} \left[\sum_{l=0}^k (-1)^l {\nu-1\choose k-l} {-\mu \choose l}
\frac{b^l}{(a-b)^l)}\right] \frac{\epsilon^{k-\rho+1}}{(k-\rho+1)}.
\end{split}
\end{equation}

The first $(\lfloor\rho\rfloor-1)$ terms of the third term of \eqref{boy} diverges as $\epsilon\rightarrow 0$. The converging and diverging parts are given by
\begin{equation}\label{conv}
\begin{split}
\tilde{C}_{\epsilon} =& \int_0^{c} \frac{x^{\nu-1}(a-x)^{-\mu}}{(b-x)^{\rho}}\mathrm{d}x \\
&+\frac{b^{\nu-1}}{(a-b)^{\mu}} \sum_{k=0}^{\infty}
\frac{(-1)^k}{b^k} \left[\sum_{l=0}^k (-1)^l {\nu-1\choose k-l} {-\mu \choose l}
\frac{b^l}{(a-b)^l)}\right] \frac{(b-c)^{k-\rho+1}}{(k-\rho+1)}\\
&-\frac{b^{\nu-1}}{(a-b)^{\mu}} \sum_{k=\lfloor\rho\rfloor}^{\infty}
\frac{(-1)^k}{b^k} \left[\sum_{l=0}^k (-1)^l {\nu-1\choose k-l} {-\mu \choose l}
\frac{b^l}{(a-b)^l)}\right] \frac{\epsilon^{k-\rho+1}}{(k-\rho+1)}
\end{split}
\end{equation} 
\begin{equation}\label{div}
\tilde{D}_{\epsilon}=-\frac{b^{\nu-1}}{(a-b)^{\mu}} \sum_{k=0}^{\lfloor\rho\rfloor}
\frac{(-1)^k}{b^k} \left[\sum_{l=0}^k (-1)^l {\nu-1\choose k-l} {-\mu \choose l}
\frac{b^l}{(a-b)^l)}\right] \frac{\epsilon^{k-\rho+1}}{(k-\rho+1)}
\end{equation}
Comparing the integral around $\mathrm{C}_{\epsilon}$ \eqref{bebe} and the diverging part \eqref{div}, we find that the contour integral assumes the form
\begin{equation}\label{bebe2}
\begin{split}
\frac{e^{i\pi(\rho-\nu)}}{2 i \sin(\pi \rho)} \int_{\mathrm{C}_{\epsilon}}
\frac{z^{\nu-1} (a+z)^{-\mu}}{(b+z)^{\rho}}\,\mathrm{d}z = -\tilde{D}_{\epsilon}+ O(\epsilon^{1-\rho+\lfloor\rho\rfloor}),
\end{split}
\end{equation} 
with the second term vanishing in the limit $\epsilon\rightarrow 0$. Then
\begin{equation}\label{meme2}
\begin{split}
& \lim_{\epsilon\rightarrow 0}\left[\int_0^{b-\epsilon} \frac{x^{\nu-1}(a-x)^{-\mu}}{(b-x)^{\rho}}\mathrm{d}x+ \frac{e^{i\pi(\rho-\nu)}}{2i \sin(\pi\rho)} \int_{C_{\epsilon}}\frac{z^{\nu-1}(a+z)^{-\mu}}{(b+z)^{\rho}} \,\mathrm{d}z\right] \\
&\hspace{30mm}= \lim_{\epsilon\rightarrow 0}\left[\int_0^{b-\epsilon} \frac{x^{\nu-1}(a-x)^{-\mu}}{(b-x)^{\rho}}\mathrm{d}x-\tilde{D}_{\epsilon} + O(\epsilon^{1-\rho+\lfloor\rho\rfloor})\right].
\end{split}
\end{equation}
Since the term $O(\epsilon^{1-\rho+\lfloor\rho\rfloor})$ vanishes in the limit, we recognize that the limit is just the finite part of the divergent integral $\int_0^b x^{\nu-1} (a+x)^{-\mu} (b-x)^{-\rho}\,\mathrm{d}x$, whose divergence arises from a non-integrable singularity at $x=b$.

Observe that the terms involving $\epsilon$ in equations \eqref{meme} and \eqref{boy} are just the negative of each other, so that the $\epsilon$ terms cancel out when \eqref{meme} and \eqref{boy} are added together. Then
\begin{equation}\label{meme2x}
\begin{split}
&\int_0^{b-\epsilon} \frac{x^{\nu-1}(a-x)^{-\mu}}{(b-x)^{\rho}}\mathrm{d}x+ \frac{e^{i\pi(\rho-\nu)}}{2i \sin(\pi\rho)} \int_{C_{\epsilon}}\frac{z^{\nu-1}(a+z)^{-\mu}}{(b+z)^{\rho}} \,\mathrm{d}z=\int_0^{c} \frac{x^{\nu-1}(a-x)^{-\mu}}{(b-x)^{\rho}}\mathrm{d}x \\
&\hspace{12mm}+\frac{b^{\nu-1}}{(a-b)^{\mu}} \sum_{k=0}^{\infty}
\frac{(-1)^k}{b^k} \left[\sum_{l=0}^k (-1)^l {\nu-1\choose k-l} {-\mu \choose l}
\frac{b^l}{(a-b)^l)}\right] \frac{(b-c)^{k-\rho+1}}{(k-\rho+1)},
\end{split}
\end{equation}
revealing that the left hand side is in fact independent of $\epsilon$. However, it now appears that the right hand side depends on $c$, which contradicts the fact that the finite part does not depend on $c$. To resolve this apparent contradiction, we take the derivative of the right hand side of \eqref{meme2x} with respect to $c$. We find that the derivative is equal to zero, meaning that the left hand side of \eqref{meme2x} is in fact independent of $c$, provided, of course, that $c$ is chosen such that the infinite series converges. Comparing the convergent part $\tilde{C}_{\epsilon}$ \eqref{conv} and the right hand side of equation \eqref{meme2x}, we find that the right hand side is just the limit of the convergent part as $\epsilon\rightarrow 0$, so that it is itself a representation of the finite part integral $\bbint{0}{b} x^{\nu-1} (a+x)^{-\mu} (b-x)^{-\rho}\,\mathrm{d}x$. However, this representation is not useful to our intentions and we now wish to evaluate the finite part explicitly without reference to the parameter $c$.  

The second term of \eqref{sese} introduces two possible cases depending on the value of $\rho$: If $0<\rho<1$, the integral is just a regular convergent integral; if $\rho\geq 1$, the integral is the finite part of a divergent integral. Nevertheless, we can simultaneously consider the two cases because the finite part coincides with the regular value of the integral when it converges. We extract the finite part by transferring the non-integrable singularity to the origin. Replacing the upper limit with $b-\epsilon$ and performing a change in variable from $x$ to $b-x$ in accordance with our discussion in Section-\ref{finitepartintegration}, we obtain the equality
\begin{equation}\label{eqfpi}
\bbint{0}{b} \frac{x^{\nu-1}(a-x)^{-\mu}}{(b-x)^{\rho}}\mathrm{d}x=\bbint{0}{b} \frac{(b-x)^{\nu-1} \, (a-b+x)^{-\mu}}{x^{\rho}} \, \mathrm{d}x.
\end{equation}
To obtain the finite part, we temporarily impose that $b/(a-b)<1$. This allows us to perform a binomial expansion of the integrand in the right hand side of \eqref{eqfpi}, and perform a term by term integration with the lower limit replaced with $\epsilon$. The result is
\begin{equation}\label{sda}
\begin{split}
&\int_{\epsilon}^{b} \frac{(b-x)^{\nu-1} \, (a-b+x)^{-\mu}}{x^{\rho}} \, \mathrm{d}x = \frac{b^{\nu-\rho}}{(a-b)^{\mu}} \sum_{k=0}^{\infty} {-\mu\choose k} \frac{1}{(a-b)^k}\\
&\hspace{22mm} \times \sum_{l=0}^{\infty} (-1)^l {\nu-1\choose l} \frac{1}{(k+l-\rho+1) b^k} \left(b^{k+l+1-\rho}-\epsilon^{k+l+1-\rho}\right)
\end{split}
\end{equation} 
The terms involving $\epsilon$ either diverge or converge to zero in the limit $\epsilon\rightarrow 0$. Those that diverge are dropped in the extraction of the finite part. Then the finite part is constituted only by the terms independent of $\epsilon$. After dropping the terms in $\epsilon$, the inner summation becomes independent of $b$ and it evaluates to
\begin{equation}\label{sume}
\sum_{l=0}^{\infty} \frac{(-1)^l}{(k+l-\rho+1)} {\nu-1\choose l} = \frac{\Gamma(\nu) \Gamma(k+2-\rho)}{(k-\rho+1) \Gamma(k+\nu-\rho+1)} .
\end{equation}
Substituting \eqref{sume} back into the finite part of \eqref{sda} and performing hypergeometric summation, the progenic finite part integral evaluates to
\begin{equation}\label{eqfpi2}
\begin{split}
\bbint{0}{b} \frac{x^{\nu-1}(a-x)^{-\mu}}{(b-x)^{\rho}}\mathrm{d}x=\frac{b^{\nu-\rho}}{(a-b)^{\mu}} \frac{\Gamma(\nu) \Gamma(1-\rho)}{\Gamma(\nu-\rho+1)} \pFq{2}{1}{\mu,1-\rho}{\nu-\rho+1}{\frac{b}{b-a}},
\end{split}
\end{equation}
for $\mathrm{Re}(\rho)>0$, $\mathrm{\nu}>0$, $\rho\notin\mathbb{Z}$, $(\rho-\nu)\notin\mathbb{Z}$ and $b/(a-b)<1$. 

Let us go back to equation \eqref{sese} and evaluate the infinite series involving finite part integrals. There are two possible cases. The first case is when $M=\lfloor\mathrm{Re}(\nu-\rho)\rfloor\geq 0$ and the second case is when $\lfloor\mathrm{Re}(\nu-\rho)\rfloor < 0$. In the former case, the integral $
\int_0^{\infty}(a+x)^{-\mu} x^{-(k+\rho-\nu+1)}\,\mathrm{d}x$
 converges for all $k=0,\dots M$ and diverges for all $k=M+1, M+2,\dots$. Under this condition, the infinite series in \eqref{sese} splits into two terms, the first is constituted by the terms $k=0,\dots,M$ at which the finite part integrals take on the convergent values of the integrals, and the rest of the integrals in the series are finite part integrals. 
 \begin{equation}
 \begin{split}\label{sese2}
 &\int_{0}^{\infty}  \frac{x^{\nu-1} \, (a+x)^{-\mu}}{(b+x)^{\rho}} \, \mathrm{d}x  =   \sum_{k=0}^{M} {- \rho \choose k} \, b^{k} \, \int_{0}^{\infty} \frac{(a+x)^{-\mu}}{x^{k+\rho-\nu+1}} \,\mathrm{d}x
 \\
 &\hspace{6mm}+\sum_{k=M+1}^{\infty} {- \rho \choose k} \, b^{k} \, \bbint{0}{\infty} \frac{(a+x)^{-\mu}}{x^{k+\rho-\nu+1}} \,\mathrm{d}x+\frac{\sin(\pi\rho)}{\sin(\pi(\rho-\nu))} \bbint{0}{b} \frac{x^{\nu-1}(a-x)^{-\mu}}{(b-x)^{\rho}}\mathrm{d}x,
 \end{split}
 \end{equation} 
 But we recall that, for a branch point singularity at the origin for the fundamental finite part, the finite part extends to the convergent case. Then the first two terms can be combined back into a single series, with the integrals given by the finite part integral \eqref{finxx} with the substitutions $\upsilon\rightarrow\mu$ and $\lambda\rightarrow (k+\rho-\nu+1)$,
\begin{align}\label{fpix1}
\begin{split}
\bbint{0}{\infty} \frac{(a+x)^{-\mu}}{x^{k+\rho-\nu+1}} \, \mathrm{d}x = \frac{(-1)^{k+1}}{a^{k+\rho-\nu+\mu}} \, \frac{\pi}{\sin((\rho-\nu) \pi)} \frac{\Gamma(k+\rho+\mu-\nu)}{ \Gamma(\mu) \Gamma(k+\rho-\nu+1)},
\end{split}
\end{align}
for all $k\in\mathbb{Z}_0^+$, $\mu\neq 0$, $(\rho-\nu)\notin\mathbb{Z}$ and $\mathrm{Re}(\rho-\nu+\mu)>0$.  

Substituting all the values of the integrals back into the Stieltjes integral \eqref{sese} and the Stieltjes integral back into the integral representation \eqref{mainlemma} of the Gauss function, we obtain the identity
\begin{equation}\label{quack1}
\begin{split}
& \pFq{2}{1}{\mu,\nu}{\mu+\rho}{1-\frac{b}{a}}  = -
\frac{\pi}{\sin{((\rho-\nu )\pi)}} \frac{\Gamma(\mu+\rho)}{\Gamma(\nu)\Gamma(\mu) \Gamma(\rho-\nu+1)}\\ &\hspace{24mm}\times \left(\frac{b}{a}\right)^{\rho-\nu} \,  \pFq{2}{1}{\rho-\nu+\mu, \rho}{\rho-\nu+1}{\frac{b}{a}}
+ \frac{\sin(\rho \pi)}{\sin{((\rho-\nu) \pi)}} \left( 1-\frac{b}{a} \right)^{-\mu} \\
&\hspace{24mm}\times \frac{\Gamma(\mu+\rho) \Gamma(1-\rho)}{\Gamma(\nu-\rho+1)\Gamma(\mu-\nu+\rho)} \, \pFq{2}{1}{\mu,1-\rho}{\nu-\rho+1}{\frac{b}{b-a}},
\end{split}
\end{equation}
for $b/(a-b)<1$, $\mathrm{Re}(\rho)>0$, $\rho\notin\mathbb{Z}$, $\nu\neq 0$, $(\rho-\nu)\notin\mathbb{Z}$, $\mu\neq 0$ and $\mathrm{Re}(\rho-\nu+\mu)>0$, where a hypergeometric summation has been performed. This result has been obtained under the condition $\lfloor\mathrm{Re}(\nu-\rho)\rfloor\geq 0$ where the integrals in the first few terms in the infinite series of \eqref{sese} are convergent. When the entire series in \eqref{sese} is composed of finite part integrals, which are also given by \eqref{fpix1}, the same expression \eqref{quack1} emmerges. Then the result \eqref{quack1} holds whether $\lfloor\mathrm{Re}(\nu-\rho)\rfloor\geq 0$ or $\lfloor\mathrm{Re}(\rho-\nu)\rfloor<0$. 

Observe that both sides of \eqref{quack1} are in terms of the variable $z=b/a$, with $0<z<1$. Replacing $b/a$ with $z$ in equation \eqref{quack1} leads to an identity for the variable $z$ which can be extended in the complex plane by letting $z$ complex. Making the replacement, we obtain the relationship
\begin{equation}\label{mainresult3}
\begin{split}
& \pFq{2}{1}{\mu,\nu}{\mu+\rho}{1-z}  = -
\frac{\pi}{\sin{((\rho-\nu )\pi)}} \frac{\Gamma(\mu+\rho)}{\Gamma(\nu)\Gamma(\mu) \Gamma(\rho-\nu+1)}\\ &\hspace{24mm}\times z^{\rho-\nu} \,  \pFq{2}{1}{\rho-\nu+\mu, \rho}{\rho-\nu+1}{z}
 + \frac{\sin(\rho \pi)}{\sin{((\rho-\nu) \pi)}} \left( 1-z \right)^{-\mu} \\
 &\hspace{24mm}\times \frac{\Gamma(\mu+\rho) \Gamma(1-\rho)}{\Gamma(\nu-\rho+1)\Gamma(\mu-\nu+\rho)} \, \pFq{2}{1}{\mu,1-\rho}{\nu-\rho+1}{\frac{z}{z-1}},
\end{split}
\end{equation}
for all $|z/(1-z)|<1$, $|\mathrm{arg}(1-z)|<\pi$ and under the same conditions as for \eqref{quack1} for the parameters. The validity of \eqref{mainresult3} can be extended to $|z|<1$ by analytic continuation. 

We wish now to express the right hand side of \eqref{mainresult3} in the variable $z$. We use the transformation equation
\begin{equation}\label{transform}
    (1-z)^{-a} \pFq{2}{1}{a,b}{c}{\frac{z}{z-1}} = \pFq{2}{1}{a,c-b}{c}{z},
\end{equation}
which is valid for $|z|<1$ by analytic continuation. We apply this transformation to the second term of equation \eqref{mainresult3}, perform the substitutions
\begin{equation}
    \Gamma(\nu-\rho+1) = -\frac{\pi}{\sin(\pi(\nu-\rho))\Gamma(\rho-\nu)}, \;\;\Gamma(\rho-\nu+1) = -\frac{\pi}{\sin(\pi(\rho-\nu))\Gamma(\nu-\rho)}\nonumber
\end{equation}
\begin{equation}
    \Gamma(1-\rho)=\frac{\pi}{\sin(\pi\rho) \Gamma(\rho)}\nonumber ,
\end{equation}
and make the shift in the variable from $z$ to $1-z$. We then set $\sigma=\mu+\rho$ and eliminate $\rho$ in favor of $\sigma$. The result is
\begin{equation}\label{mainresult4x}
\begin{split}
& \pFq{2}{1}{\mu,\nu}{\sigma}{z}  
	 = \frac{\Gamma(\sigma)) \Gamma(\mu+\nu-\sigma)}{\Gamma(\nu) \Gamma(\mu)} \, (1-z)^{\sigma-\mu-\nu} \pFq{2}{1}{\sigma-\nu,\sigma-\mu}{\sigma-\mu-\nu+1}{1-z}
	\\
& \hspace{24mm}+ \frac{\Gamma(\sigma) \Gamma(\sigma-\mu-\nu)}{\Gamma(\sigma-\nu)\Gamma(\sigma-\mu)} \,
\pFq{2}{1}{\mu,\nu}{\mu+\nu-\rho+1}{1-z},
\end{split}
\end{equation}
which is valid for $|1-z|<1$, $|\mathrm{Re}(1-z)|<\pi$, $\mathrm{Re}(\sigma-\mu)>0$, $(\sigma-\mu)\notin\mathbb{Z}$, $(\sigma-\mu-\nu)\notin\mathbb{Z}$ and $\mathrm{Re}(\sigma-\nu)>0$. The restrictions on the parameters may be lifted by analytic continuation provided $(\sigma-\mu-\nu)\notin \mathbb{Z}$. Equation \eqref{mainresult4x} is a known transformation equation for the Gauss function \cite{transform1}.

\subsubsection{Case $\rho\in\mathbb{Z}^+$}
We obtained equation \eqref{mainresult3} under the condition that $\rho\notin\mathbb{Z}^+$ so that the transformation equation \eqref{mainresult4x} was established only (through \eqref{mainresult3}) for $\rho=(\sigma-\mu)\notin \mathbb{Z}^+$. We now establish that \eqref{mainresult4x} holds even for positive integer $\rho$. Since $(\rho-\nu)\notin\mathbb{Z}$ and $\rho\in\mathbb{Z}^+$, it is now necessary that $\nu\notin\mathbb{Z}$. Let $\rho=n\in\mathbb{Z}^+$.  This time the kernel has a pole of order $n$ at $z=-b$. Since the singularity at the origin is a branch point singularity due to the non-integral value of $\nu$, the Stieltjes integral is extracted from
\begin{equation}
\frac{1}{e^{2\pi i\nu}-1} \int_{\mathrm{C}} \frac{z^{\nu-1} (a+z)^{-\mu}}{(b+z)^{\rho}} \, \mathrm{d}z,
\end{equation}
where the pre-factor follows from $(e^{-2\pi i (n-\nu)}-1)=(e^{2\pi i \nu}-1)$. Deforming the contour $\mathrm{C}$ to $\mathrm{C}'$ as shown in Figure-\ref{deformation} (bottom figure), we obtain
\begin{equation}
\begin{split}
\int_{0}^{\infty} \frac{x^{\nu-1} \, (a+x)^{-\mu}}{(b+x)^{n}} \,\mathrm{d}x
=& \frac{1}{e^{2\pi i \nu}-1} \int_C \frac{z^{\nu-1}(a+z)^{-\mu}}{(b+z)^{n}} \,\mathrm{d}z
\\
&\hspace{4mm}- \frac{2 \pi i}{e^{2\pi i \nu}-1} \, \mathrm{Res} \left[ \frac{z^{\nu-1} \, (a+z)^{-\mu}}{(b+z)^{n}}\right]_{z=-b}.
\end{split}
\end{equation}
The residue term arises from the fact that $z=-b$ is now a pole of the kernel in contrast to the previous where it is a branch point. Here the singular contribution arises from a pole singularity. Again expanding the kernel in the first term and distributing the contour integration, the resulting integrals are again identified to be finite part integrals and we arrive at the evaluation,
\begin{align}\label{representation1a}
\begin{split}
\int_{0}^{\infty} \frac{x^{\nu-1} \, (a+x)^{-\mu}}{(b+x)^{n}} \,\mathrm{d}x
= & \sum_{k=0}^{\infty} {- n \choose k} \, b^{k} \,\bbint{0}{\infty} \frac{(a+x)^{-\mu}}{x^{k+n-\nu+1}} \,\mathrm{d}x \\
&\hspace{12mm}
- \frac{2 \pi i}{e^{2\pi i \nu}-1} \, \mathrm{Res} \left[ \frac{z^{\nu-1} \, (a+z)^{-\mu}}{(b+z)^{n}}\right]_{z=-b},
\end{split}
\end{align}
which is valid for $\mathrm{Re}(\nu)>0$, $\nu\notin\mathbb{Z}$, $\mu\neq 0$, $n\in\mathbb{Z}^+$, $\mathrm{Re}(n+\mu-\nu)>0$ and $a>b>0$.

We have the same situation as in the previous case that we have two separate cases: the first spliting the infinite series in two invoving convergent integrals and the rest finite part integrals, and the second case involves only finite part integrals. We likewise find that the finite part integral extends even in the case of convergent integrals, so that we end up with just one series. The relevant finite part integral in this case is given by \eqref{finxx} with the substitutions $\upsilon\rightarrow\mu$ and $\lambda\rightarrow (k+n-\nu+1)$, 
\begin{equation}\label{finitepart1}
\begin{split}
\bbint{0}{\infty} \frac{(a+x)^{-\mu}}{ x^{k+n-\nu+1}} \,\mathrm{d}x
= \frac{(-1)^{k+n}\pi \,\Gamma(k+n+\mu-\nu)}{a^{k+n+\mu-\nu}\sin(\pi\nu)\,\Gamma(\mu)\,\Gamma(k+\mu-\nu+1)}
\end{split}
\end{equation}
for all $k\in\mathbb{Z}^+_0$ under the conditions for \eqref{representation1a}. On the other hand, the residue is obtained from the known limit representation for the residue,
\begin{equation}\label{res}
\begin{split}
\mathrm{Res}  \left[ \frac{z^{\nu-1} \, (a+z)^{-\mu}}{(b+z)^{n}}\right]_{z=-b} = \frac{1}{(n-1)!}\lim_{z\rightarrow -b} \frac{\mathrm{d}^{n-1}}{\mathrm{d}z^{n-1}}\left[z^{\nu-1} \, (a+z)^{-\mu}\right] ,
\end{split}
\end{equation}
where the relevant derivatives are given by
\begin{equation}
\frac{\mathrm{d}^n}{\mathrm{d}z^n} z^{\nu-1}= (\nu-n)_n \, z^{\nu-n-1},
\end{equation}
\begin{equation}
\frac{\mathrm{d}^n}{\mathrm{d}z^n}\left(z+a\right)^{-\mu}= (-1)^n (\mu)_n (z+a)^{-\mu-n}.
\end{equation}
Then by Leibnitz rule of differentiation, the residue is calculated to be
\begin{equation}\label{res2xx}
\begin{split}
\mathrm{Res}  \left[ \frac{z^{\nu-1} \, (a+z)^{-\mu}}{(b+z)^{n}}\right]_{z=-b} = \frac{e^{i\pi \nu}}{\Gamma(n)} & \sum_{k=0}^{n -1}  {n -1 \choose k} \frac{(\nu-k)_{k} (\mu)_{n-1-k}}{b^{k-\nu+1}}
(a-b)^{k-\mu-n+1}.
\end{split}
\end{equation}

Collecting all the terms together and combining the first two terms and substituting the overall result back into equation \eqref{mainlemma} lead to an identity for $_2F_1$ in power series of the variable $b/a$. The resulting identity holds for $0<b/a<1$ but it admits a complex extension by replacing $b/a$ with the complex variable $z$ for $|z|<1$, which can be further analytically continued in the rest of the complex plane by hypergeometric summation. Performing the extension yields the identity   
\begin{equation}\label{mainresult1}
\begin{split}
& \pFq{2}{1}{\mu,\nu}{\mu+n}{1-z} 
= \frac{(-1)^{n} \pi \Gamma(\mu+n)\, z^{n-\nu}}{\sin(\pi \nu) \Gamma(\nu)\Gamma(\mu) \Gamma(n-\nu+1)}\pFq{2}{1}{n-\nu+\mu, \rho}{n-\nu+1}{z}
	\\
&\hspace{2mm} - \frac{(-1)^{n} \pi \Gamma(\mu+n)}{\sin(\pi \nu)\Gamma(\nu)\Gamma(\mu-\nu+n)\Gamma(n)} \sum_{k=0}^{n-1} {n-1 \choose k} (\mu)_k (\nu-n+1+k)_{n-1-k} \frac{z^k}{(1-z)^{k+\mu}} ,
\end{split}
\end{equation}
for $\mathrm{Re}(\nu)>0$, $\nu\notin\mathbb{Z}$, $\mu\neq 0$, $n\in\mathbb{Z}^+$, $\mathrm{Re}(n+\mu-\nu)>0$, with all functions involved take their principal values. For a given fixed $n$, the restrictions on the the parameters may be lifted by analytic continuation. 

We wish to rewrite \eqref{mainresult1} as well in powers of $z$ alone. The first term already is in powers of $z$. To rewrite the second term, we perform the following substitutions
\begin{equation*}
(\nu-n+1+k)_{n-1-k} = \frac{\Gamma(\nu)}{\Gamma(\nu-n+1) (\nu-n+1)_k},\;\;
{n-1\choose k} = (-1)^k \frac{(1-n)_k}{k!}
\end{equation*}
and we readily recognize that the sum can be evaluated in terms of the Gauss function. The result is
\begin{equation}\label{xx12}
\begin{split}
&\sum_{k=0}^{n-1} {n-1 \choose k} (\mu)_k (\nu-n+1+k)_{n-1-k} \frac{z^k}{(1-z)^{k+\mu}}\\ 
&\hspace{24mm}= \frac{\Gamma(\nu)}{\Gamma(\nu-n+1)} (1-z)^{-\mu} \pFq{2}{1}{\mu,1-n}{\nu-n+1}{\frac{z}{z-1}}
\end{split}
\end{equation}
Applying the transformation equation \eqref{transform} in the right hand side of \eqref{xxx1}, we obtain the equality 
\begin{equation}\label{xxx12}
\begin{split}
&\sum_{k=0}^{n-1} {n-1 \choose k} (\mu)_k (\nu-n+1+k)_{n-1-k} \frac{z^k}{(1-z)^{k+\mu}}=\frac{\Gamma(\nu)}{\Gamma(\nu-n+1)} \pFq{2}{1}{\mu,\nu}{\nu-n+1}{z} 
\end{split}
\end{equation}
We substitute \eqref{xxx12} back into \eqref{mainresult1} and make the following substitutions,
\begin{equation}\nonumber
\Gamma(n-\nu+1)= (-1)^n \frac{\pi}{\sin(\pi \nu) \Gamma(n-\nu)},\;\; \Gamma(\nu-n+1)= (-1)^{n+1} \frac{\pi}{\sin(\pi \nu) \Gamma(\nu-n)}.
\end{equation}
This is followed by making the shift from $z$ to $(1-z)$ and setting $\sigma=\mu+n$. This leads to the transformation \eqref{mainresult4x} with $\rho=n$ for non-integer $\nu$. 

This example demonstrates the general situation that the case for a pole at $-b$ is a special value of the analytic continuation with respect to $\rho$ of the result for the case of branch point. It is for this reason that we will restrict our consideration for the case of a branch point only at $-b$ for the rest of the paper. 

\subsection{Case $(\nu-\rho)\in\mathbb{Z}$ and $\nu,\rho\notin \mathbb{Z}$} \label{SS4.2}
Under the condition that the difference  $(\nu-\rho)$ is an integer, the divergence is due to a pole singularity at $z=0$ , but since $\rho$ is not an integer $z=-b$ is a branch point.  Since the origin is a pole, we extract the Stieltjes integral \eqref{si1} from the contour integral
\begin{equation}
\frac{1}{2\pi i}\int_{\mathrm{C}} \frac{z^{\nu-1} (a+z)^{-\mu}}{(b+z)^{\rho}}\,\left(\log z - i\pi\right)\, \mathrm{d}z ,
\end{equation}
where the contour $\mathrm{C}$ is the same as above. Again collapsing the contour into the real line through the contour $\mathrm{C}'$ yields the representation for the integral
\begin{equation}
\begin{split}
&\int_0^{\infty} \frac{x^{\nu-1} (a+x)^{-\mu}}{(b+x)^{\rho}} \mathrm{d}x=  \frac{1}{2\pi i} \int_C \frac{z^{\nu-1} (a+z)^{-\mu}}{(b+z)^{\rho}}\,\left(\log z - i\pi\right)\, \mathrm{d}z \\
& \hspace{24mm}+ (-1)^{\nu-\rho+1} \frac{\sin(\pi \rho)}{\pi} \lim_{\epsilon\rightarrow 0} \left[\int_0^{b-\epsilon} \frac{x^{\nu-1} (a-x)^{-\mu}}{(b-x)^{\rho}}\,\ln(x)\,\mathrm{d}x \right. \\
& \hspace{24mm}\left.  + \frac{e^{-i\pi\nu}}{(1-e^{-2\pi\rho i})}\int_{C_{\epsilon}} \frac{z^{\nu-1} (a+z)^{-\mu}}{(b+z)^{\rho}}\,\left(\log z - i\pi\right)\, \mathrm{d}z\right],
\end{split}
\end{equation}
where the fact that $(\nu-\rho)$ is an integer has been used to arrive at the second term. Using the same method of proof applied earlier, it can be established that the limit is the finite part of the divergent integral $\int_0^b x^{\nu-1} (a+x)^{-\mu} (b+x)^{-\rho}\,\mathrm{d}x$. We perform a similar expansion in the first term and perform term by term integration to obtain the desired evaluation of the integral,
\begin{equation}
\begin{split}\label{may}
&\int_0^{\infty} \frac{x^{\nu-1} (a+x)^{-\mu}}{(b+x)^{\rho}}\mathrm{d}x =  \sum_{k=0}^{\infty} {- \rho \choose k} \, b^{k} \, \bbint{0}{\infty} \frac{(a+x)^{-\mu}}{x^{k+\rho-\nu+1}} \,\mathrm{d}x \\
& \hspace{24mm}+  (-1)^{\nu-\rho+1}\frac{\sin(\pi\rho)}{\pi} \bbint{0}{b} \frac{x^{\nu-1} \, (a-x)^{-\mu}}{(b-x)^{\rho}} \, \ln(x)\, \mathrm{d}x .
\end{split}
\end{equation}
which is valid for $a>b>0$, $\mathrm{Re}(\nu)>0$, $\mu\neq 0$, $(\rho+\mu-\nu)>0$, $\mathrm{Re}(\rho)>0$, $\rho\notin\mathbb{Z}$ and $(\nu-\rho)\in\mathbb{Z}$. The singular contribution in this case is again a progenic finite-part integral. 

Let us now obtain the progenic finite part integral
(or the regular value when the integral converges) in equation \eqref{may}.  Let $\epsilon$ be some positive arbitrarilly small number and consider the integral
\begin{equation}
\frac{a^{\mu}}{b^{\nu-\rho}}\int_{0}^{b-b\epsilon} \frac{x^{\nu-1} (a-x)^{-\mu}}{(b-x)^{\rho}}\,\ln(x)\, \mathrm{d}x = C_{\epsilon} + D_{\epsilon},
\end{equation}
where $C_{\epsilon}$ and $D_{\epsilon}$ are the convergent and divergent parts, respectively, in the limit as $\epsilon$ approaches zero. Again the desired finite part is the limit $\lim_{\epsilon\rightarrow 0} C_{\epsilon}$. We change variable $x=b y$. Under the condition that $b>0$, the integral takes on the form
\begin{equation}\label{integralorig}
\begin{split}
&\frac{a^{\mu}}{b^{\nu-\rho}}\int_0^{b-b\epsilon} \frac{x^{\nu-1} (a-x)^{-\mu}}{(b-x)^{\rho}}\,\ln(x)\, \mathrm{d}x =  \ln(b) \int_0^{1-\epsilon} \frac{y^{\nu-1} (1-\frac{b}{a} y)^{-\mu}}{(1-y)^{\rho}}\, \mathrm{d}y \\
& \hspace{20mm}+\int_0^{1-\epsilon} \frac{y^{\nu-1} \left(1-\frac{b}{a} y\right)^{-\mu}}{(1-y)^{\rho}}\,\ln(y)\,\mathrm{d}y .
\end{split}
\end{equation}
This expression implies that the desired finite part decomposes into a combination of two finite part integrals,
\begin{equation}\label{integralorigx}
\begin{split}
&\frac{a^{\mu}}{b^{\nu-\rho}}\bbint{0}{b} \frac{x^{\nu-1} (a-x)^{-\mu}}{(b-x)^{\rho}}\,\ln(x)\, \mathrm{d}x =  \ln(b) \bbint{0}{1} \frac{y^{\nu-1} (1-\frac{b}{a} y)^{-\mu}}{(1-y)^{\rho}}\, \mathrm{d}y \\
& \hspace{20mm}+\bbint{0}{1} \frac{y^{\nu-1} \left(1-\frac{b}{a} y\right)^{-\mu}}{(1-y)^{\rho}}\,\ln(y)\,\mathrm{d}y .
\end{split}
\end{equation}

To obtain the finite part integrals in the right hand side of equation \eqref{integralorigx}, we need to evaluate the integrals 
\begin{equation}
\int_0^{1-\epsilon} \frac{y^{\nu-1} (1-\frac{b}{a} y)^{-\mu}}{(1-y)^{\rho}}\,g_l(y) \mathrm{d}y = C_{\epsilon}^{(l)} + D_{\epsilon}^{(l)}, \;\;\; l=1, 2, 
\end{equation}
where $g_1(y)=1$ and $g_2(y)=\ln(y)$, and the  $C_{\epsilon}^{(l)}$'s and $D_{\epsilon}^{(l)}$'s are the respective converging and diverging parts of the integrals. To evaluate these integrals, we let $b/a<(1-\epsilon)$ and expand binomially the factor $(1-b y/a)^{-\mu}$, followed by term by term integration, which is allowed because the series converges uniformly along the interval of integration $[0,1-\epsilon]$. Then 
\begin{equation}
\int_0^{1-\epsilon} \frac{y^{\nu-1} (1-\frac{b}{a} y)^{-\mu}}{(1-y)^{\rho}}\,g_l(y)\, \mathrm{d}y = \sum_{k=0}^{\infty} (-1)^k {-\mu \choose k}\left(\frac{b}{a}\right)^k \int_0^{1-\epsilon} \frac{y^{\nu-1+k}}{(1-y)^{\rho}}\,g_l(y)\,\mathrm{d}y .
\end{equation}
We are lead to evaluate the integrals
\begin{equation}\label{rawintegral}
\int_0^{1-\epsilon} \frac{y^{\nu+k-1}}{(1-y)^{\rho}}\,g_l(y)\,\mathrm{d}y=\tilde{C}^{(l)}_{k,\epsilon} + \tilde{D}^{(l)}_{k,\epsilon},
\end{equation}
where the $\tilde{C}^{(l)}_{\epsilon}$'s and $\tilde{D}^{(l)}_{\epsilon}$'s are the converging and diverging parts of the integrals for a given $k$. 
Clearly the converging part $C_{\epsilon}^{(1)}$ accumulates from the converging parts $\tilde{C}_{k,\epsilon}$'s. These imply that the desired finite integrals are given by
\begin{equation}\label{fpixxx}
\bbint{0}{1}\frac{y^{\nu-1} (1-\frac{b}{a} y)^{-\mu}}{(1-y)^{\rho}}\, \mathrm{d}y = \sum_{k=0}^{\infty} (-1)^k {-\mu \choose k}\left(\frac{b}{a}\right)^k \bbint{0}{1} \frac{y^{\nu-1+k}}{(1-y)^{\rho}}\,\mathrm{d}y ,
\end{equation}
\begin{equation}\label{fpixxx2}
\bbint{0}{1}\frac{y^{\nu-1} (1-\frac{b}{a} y)^{-\mu}}{(1-y)^{\rho}}\,\ln(y)\, \mathrm{d}y = \sum_{k=0}^{\infty} (-1)^k {-\mu \choose k}\left(\frac{b}{a}\right)^k \bbint{0}{1} \frac{y^{\nu-1+k}}{(1-y)^{\rho}}\,\ln(y)\,\mathrm{d}y .
\end{equation}
The two infinite series are convergent. We are finally led to evaluate the finite part integrals
\begin{equation}\label{toevaluate}
\bbint{0}{1} \frac{y^{\nu-1+k}}{(1-y)^{\rho}}\,\mathrm{d}y,\;\;\;
\bbint{0}{1} \frac{y^{\nu-1+k}}{(1-y)^{\rho}}\,\ln(y)\,\mathrm{d}y .
\end{equation}

Let us evaluate the first of \eqref{toevaluate}. For $\sigma$ and $\rho$ with $\mathrm{Re}(\sigma)>0$, $\mathrm{Re}(\rho)>0$ and $\rho\notin\mathbb{Z}^+$, we consider the finite part integral
\begin{equation}
\bbint{0}{1}\frac{y^{\sigma-1}}{(1-y)^{\rho}}\,\mathrm{d}y .
\end{equation} Let $0<\epsilon<1$ and consider the integral
\begin{equation}
\int_0^{1-\epsilon} \frac{y^{\sigma-1}}{(1-y)^{\rho}}\,\mathrm{d}y .
\end{equation}  To evaluate this integral, we use two facts about the incomplete beta function, its integral representation
\begin{equation}\label{betaintegrep}
B_{z}(a,b) = \int_0^z t^{a-1} (1-t)^{b-1}\,\mathrm{d}t,\; \mathrm{Re}(a)>0,
\end{equation}
and its behavior in the neighbhood of $z=1$,
\begin{equation}\label{betainc1}
B_z(a,b)= B(a,b) - \frac{(1-z)^{b} z^a}{b}\left(1+O(z-1)\right),
\end{equation}
where $B(a,b)$ is the beta function. Comparing the integral \eqref{rawintegral} and the integral representation of the incomplete beta function \eqref{betaintegrep}, we find that integral \eqref{rawintegral} evaluates to the incomplete beta function with $z=1-\epsilon$, $a=\sigma$ and $b=1-\rho$. Then applying the expansion \eqref{betainc1} to the result yields
\begin{equation}
\int_0^{1-\epsilon} \frac{y^{\sigma-1}}{(1-y)^{\rho}}\,\mathrm{d}y = \frac{\pi\Gamma(\sigma)}{\sin(\pi\rho) \Gamma(\rho) \Gamma(1+\sigma-\rho)} + O(\epsilon^{1-\rho}),
\end{equation}
where we have made a simplification under the condition that $\rho$ is a non-integer. If $\mathrm{Re}(\rho)<1$, the second term vanishes as $\epsilon\rightarrow 0$ and we have the case of a converging integral, and the finite part is just the value of the convergent integral. If, on the other hand, $\mathrm{Re}(\rho)>1$, the second term diverges and the finite part is obtained by dropping the diverging part and taking the limit of the converging part. Thus we arrive at the finite part integral
\begin{equation}\label{fpingteg1x}
\bbint{0}{1}\frac{y^{\sigma-1}}{(1-y)^{\rho}}\,\mathrm{d}y = \frac{\pi\Gamma(\sigma)}{\sin(\pi\rho) \Gamma(\rho) \Gamma(1+\sigma-\rho)},
\end{equation}
which is valid for $\mathrm{Re}(\sigma)>0$, $\mathrm{Re}(\rho)>0$ and $\rho\notin\mathbb{Z}^+$. Since $1/\Gamma(z)=0$ for $z=0,-1, -2, \dots$, the finite part assumes the zero value for $1+\sigma-\rho=0,-1,-2,\dots$. Explicitly, we have
\begin{equation}\label{fpingteg1xxx}
\bbint{0}{1}\frac{y^{\sigma-1}}{(1-y)^{\rho}}\,\mathrm{d}y = 0,
\end{equation}
for $(\sigma-\rho)\in\mathbb{Z}^-$, $\mathrm{Re}(\sigma)>0$, $\mathrm{Re}(\rho)>0$ and $\rho\notin\mathbb{Z}^+$.

Now observe that the integral $\int_0^{1-\epsilon}y^{\sigma-1}\ln(y) (1-y)^{-\rho}\mathrm{d}y$  is just the derivative of the integral $\int_0^{1-\epsilon}y^{\sigma-1} (1-y)^{-\rho}\mathrm{d}y$ with respect to the parameter $\sigma$, which follows from the fact that the integral can be differentiated under the integral with respect to $\sigma$ because the later integral converges uniformly in the entire integration interval $[0,1-\epsilon]$. This implies that the second finite part in \eqref{toevaluate} can be obtained from \eqref{fpingteg1x} by differentiation with respect to $\sigma$. Performing the indicated differentiation, we obtain the finite part integral
\begin{equation}\label{fpingteg1xxx2}
\bbint{0}{1}\frac{y^{\sigma-1}}{(1-y)^{\rho}}\,\ln(y)\,\mathrm{d}y = \frac{\pi\Gamma(\sigma)}{\sin(\pi\rho) \Gamma(\rho) \Gamma(1+\sigma-\rho)} \left(\psi(\sigma)-\psi(\sigma-\rho+1)\right),
\end{equation}
for $\mathrm{Re}(\sigma)>0$, $\mathrm{Re}(\rho)>0$ and $\rho\notin\mathbb{Z}^+$. For $\sigma-\rho+1=0,-1,-2,\dots$, we have a removable singularity. Using the fact that  $\lim_{z\rightarrow -n} \psi(z)/\Gamma(z) = (-1)^{n+1} n!$ for a positive integer $n$, we have the following special values for the finite part integral
\begin{equation}\label{fpingteg1xx3}
\bbint{0}{1}\frac{y^{\sigma-1}}{(1-y)^{\rho}}\,\ln(y)\,\mathrm{d}y = (-1)^{\rho-\sigma+1} \frac{\pi\Gamma(\sigma) (\rho-\sigma-1)!}{\sin(\pi\rho) \Gamma(\rho)},
\end{equation}
for $(\sigma-\rho)\in\mathbb{Z}^-$, $\mathrm{Re}(\sigma)>0$, $\mathrm{Re}(\rho)$ and $\rho\notin\mathbb{Z}^+$.

From the finite part integrals \eqref{fpingteg1x}, \eqref{fpingteg1xxx}, \eqref{fpingteg1xxx2}, \eqref{fpingteg1xx3}, it is clear that the explicit form of the desired finite part integral,
\begin{equation}
\bbint{0}{b} \frac{x^{\nu-1} \, (a-x)^{-\mu}}{(b-x)^{\rho}} \, \ln(x)\, \mathrm{d}x,
\end{equation}
depends on the relative difference of the parameters $\rho$ and $\nu$. We now consider separately the cases $(\nu-\rho)\in\mathbb{Z}^+$ and $(\rho-\nu)\in\mathbb{Z}^+_0$.  

\subsubsection{Case $(\nu-\rho)\in\mathbb{Z}^+$} Under this case, the expansion decomposes into a finite series of convergent integrals and an infinite series of finite part integrals. We have the expansion
\begin{equation}
\begin{split}
&\int_0^{\infty} \frac{x^{\nu-1} (a+x)^{-\mu}}{(b+x)^{\rho}} \mathrm{d}x  =  \sum_{k=0}^{(\nu-\rho)-1} \, {- \rho \choose k} \, b^{k} \int_0^{\infty} \ \frac{(a+x)^{-\mu}}{x^{k+\rho-\nu+1}} \,\mathrm{d}x \\
& + \sum_{k=(\nu-\rho)}^{\infty} {- \rho \choose k} \, b^{k} \, \bbint{0}{\infty} \frac{(a+x)^{-\mu}}{x^{k+\rho-\nu+1}} \,\mathrm{d}x
- (-1)^{\nu-\rho}\frac{\sin(\pi\rho)}{\pi} \bbint{0}{b} \frac{x^{\nu-1} \, (a-x)^{-\mu}}{(b-x)^{\rho}} \, \ln(x)\, \mathrm{d}x,
\label{case4b}
\end{split}
\end{equation}
which is valid for $a>b>0$, $\mathrm{Re}(\nu)>0$, $\mu\neq 0$, $(\rho+\mu-\nu)>0$, $\mathrm{Re}(\rho)>0$, $\rho\notin\mathbb{Z}$ and $(\nu-\rho)\in\mathbb{Z}$.  

First, let us simplify the first two group of terms of equation \eqref{case4b}. The integrals in the first group of terms evaluate are given by \eqref{finiteres} with the substitutions $\upsilon\rightarrow\mu$ and $\lambda\rightarrow (k+\rho-\nu+1)$,
\begin{equation}\label{kwe}
    \int_0^{\infty} \frac{(a+x)^{-\mu}}{x^{k+\rho-\nu+1}}\,\mathrm{d}x = \frac{\Gamma(\nu-\rho-k)\Gamma(k+\mu-\nu+\rho)}{a^{k+\mu-\nu+\rho} \Gamma(\mu)},
\end{equation}
for $k=0, \dots, (\nu-\rho-1)$, $\mu\neq 0$ and $\mathrm{Re}(\rho-\nu+\mu)>0$. Substituting the values of the convergent integrals, the first term evaluates to
\begin{equation}\label{xxx1}
\begin{split}
    &\sum_{k=0}^{(\nu-\rho)-1} \, {- \rho \choose k} \, b^{k} \int_0^{\infty} \ \frac{(a+x)^{-\mu}}{x^{k+\rho-\nu+1}} \,\mathrm{d}x \\
    &\hspace{8mm}= \frac{\Gamma(\nu-\rho) \Gamma(\mu-\nu+\rho)}{a^{\mu+\rho-\nu}\Gamma(\mu)} \sum_{k=0}^{\nu-\rho-1} \frac{(\rho)_k (\mu-\nu+\rho)_k}{(1-\nu+\rho)_k k!} \left(\frac{b}{a}\right)^k,
       \end{split}
\end{equation}
under the conditions for \eqref{kwe}, in which the identity $\Gamma(z-n)= (-1)^n \Gamma(z)/(1-z)_n$ for positive integer $n$ has been applied to arrive at the expression.

On the other hand, the finite part integrals in the second term are given by equation \eqref{FPIpole} with the substitutions $\upsilon\rightarrow\mu$, $s\rightarrow a$, $n\rightarrow (k+\rho-\nu)$. The result is
\begin{equation}\label{firstterm}
\begin{split}
\bbint{0}{\infty} \frac{(a+x)^{-\mu}}{x^{k+\rho-\nu+1}} \,\mathrm{d}x &= \frac{(-1)^{k+\rho-\nu}}{a^{k+\rho-\nu+\mu}} \frac{(\mu)_{k+\rho-\nu}}{(k+\rho-\nu)!} \\
&\hspace{10mm}\times\left[\ln(a) + \psi(k+\rho-\nu+1)-\psi(k+\rho-\nu+\mu)\right],
\end{split}
\end{equation}
for all $k=(\nu-\rho), (\nu-\rho+1),\dots$, $\mu\neq 0$ and $\mathrm{Re}(\rho-\nu+\mu)>0$. For the second group of terms, we substitute the finite part integrals and shift the index of summation from $k$ to $k+\nu-\rho$. The coefficient of $\ln(a)$ is readily evaluated in terms of a hypergeometric function. Then, after writing all binomial coefficients in terms of the gamma function, the second term assumes the form
\begin{equation}
    \begin{split}
        &\sum_{k=(\nu-\rho)}^{\infty} {- \rho \choose k} \, b^{k} \, \bbint{0}{\infty} \frac{(a+x)^{-\mu}}{x^{k+\rho-\nu+1}} \,\mathrm{d}x\\
        &\hspace{4mm}= \frac{b^{\nu-\rho}}{a^{\mu}}(-1)^{\nu-\rho} \frac{\Gamma(\nu)}{\Gamma(\rho) \Gamma(\nu-\rho+1)} \pFq{2}{1}{\mu,\nu}{\mu-\rho}{\frac{b}{a}} \,\ln(a) \\
        &\hspace{6mm}+\frac{b^{\nu-\rho}}{a^{\mu}} \frac{(-1)^{\nu-\rho}}{\Gamma(\mu) \Gamma(\rho)} \sum_{k=0}^{\infty} \frac{\Gamma(k+\mu) \Gamma(k+\nu)}{\Gamma(k+\nu-\rho+1) k!} \left[\psi(k+1)-\psi(k+\mu)\right]\left(\frac{b}{a}\right)^k ,
    \end{split}
\end{equation}
which is valid for  $b/a<1$, $(\nu-\rho)\in\mathbb{Z}^+$ and $\mathrm{Re}(\rho-\nu+\mu)>0$.

We now evaluate the progenic finite part integral in equation \eqref{case4b}. Under the condition that $(\nu-\rho)$ is a positive integer and with $\sigma=\nu+k$, the finite part integral \eqref{fpingteg1x} does not vanish for all terms in the expansion, so that
\begin{equation}\label{be}
\bbint{0}{1}\frac{y^{k+\nu-1}}{(1-y)^{\rho}}\,\mathrm{d} y = \frac{\pi \Gamma(k+\nu)}{\sin(\pi\rho) \Gamma(\rho) \Gamma(k+\nu-\rho+1)},
\end{equation} 
for all $k\in\mathbb{Z}^+_0$, $\mathrm{Re}(\nu)>0$, $\mathrm{Re}(\rho)>0$ and $\rho\notin\mathbb{Z}^+$. Substituting these finite parts \eqref{be} back into equation \eqref{fpixxx} and writing the gamma functions in terms of the pochhamer symbol yield  the finite part integral 
\begin{equation}\label{integralseries}
    \bbint{0}{1} \frac{y^{\nu-1} (1-\frac{b}{a} y)^{-\mu}}{(1-y)^{\rho}}\, \mathrm{d}y =\frac{\pi \Gamma(\nu)}{\sin(\pi\rho)\Gamma(\rho) \Gamma(\nu-\rho+1)}\pFq{2}{1}{\mu,\nu}{\nu-\rho+1}{\frac{b}{a}},
\end{equation}
for $b/a<1$ and under the conditions for \eqref{be}. 

Similarly the logarithmic finite part integral takes uniform expression for all $k$, 
\begin{equation}\label{fpingteg1xx}
\bbint{0}{1}\frac{y^{\nu+k-1}}{(1-y)^{\rho}}\,\ln(y)\,\mathrm{d}y = \frac{\pi\Gamma(\nu+k)}{\sin(\pi\rho) \Gamma(\rho) \Gamma(1+k+\nu-\rho)} \left(\psi(k+\nu)-\psi(1+k+\nu-\rho)\right),
\end{equation}
for all $k\in\mathbb{Z}^+_0$, $\mathrm{Re}(\nu)>0$, $\mathrm{Re}(\rho)>0$ and $\rho\notin\mathbb{Z}^+$. Then substituting the finite part into equation \eqref{fpixxx2}, the logarithmic finite part integral is given by
\begin{equation}\label{integrallogseries}
\begin{split}
    &\bbint{0}{1} \frac{y^{\nu-1} (1-\frac{b}{a} y)^{-\mu}}{(1-y)^{\rho}}\,\ln(y)\, \mathrm{d}y=\frac{\pi}{\sin(\pi\rho) \Gamma(\mu)\Gamma(\rho)} \\ 
    &\hspace{10mm}\times\sum_{k=0}^{\infty} \frac{\Gamma(k+\mu) \Gamma(k+\nu)}{\Gamma(1+k+\nu-\rho) k!}  \left[\psi(k+\nu)-\psi(1+k+\nu-\rho)\right] \left(\frac{b}{a}\right)^k ,
    \end{split}
\end{equation}
which is valid for all $a>b>0$ and under the conditions for \eqref{fpingteg1xx}. 

We substitute the finite part integrals \eqref{integralseries} and \eqref{integrallogseries} back into equation \eqref{integralorigx}, the desired finite part is given by
\begin{equation}\label{coco}
\begin{split}
&\frac{a^{\mu}}{b^{\nu-\rho}} \bbint{0}{b} \frac{x^{\nu-1} \, (a-x)^{-\mu}}{(b-x)^{\rho}} \, \ln(x)\, \mathrm{d}x
= \frac{\pi \Gamma(\nu)\ln(b)}{\sin(\pi \rho) \Gamma(\rho) \Gamma(1+\nu-\rho)} \pFq{2}{1}{\mu,\nu}{\nu-\rho+1}{\frac{b}{a}} \\
&+ \frac{\pi}{\sin(\pi\rho) \Gamma(\mu)\Gamma(\rho)}  
\sum_{k=0}^{\infty} \frac{\Gamma(k+\mu) \Gamma(k+\nu)}{\Gamma(1+k+\nu-\rho) k!}
\left[\psi(k+\nu)-\psi(1+k+\nu-\rho)\right] \left(\frac{b}{a}\right)^k,
\end{split}
\end{equation} 
for $a>b>0$, $\mathrm{Re}(\nu)>0$, $\mu\neq 0$, $\mathrm{Re}(\rho)>0$ and $\rho\notin\mathbb{Z}$. We substitute this back into expression \eqref{case4b}. The $\ln(a)$ and $\ln(b)$ combine and so with the terms involving the digamma functions. As in the previous cases, the resulting expression is in terms of the variable $b/a$ which we now replace complex $z$ by the principle of analytic continuation. We obtain the following result,
\begin{equation}\label{repcase4bx}
    \begin{split}
        &\pFq{2}{1}{\mu,\nu}{\mu+\rho}{1-z}= \frac{\Gamma(\mu+\rho)\Gamma(\nu-\rho)}{\Gamma(\mu)\Gamma(\nu) z^{\nu-\rho}} \sum_{K=0}^{\nu-\rho-1} \frac{(\rho)_k (\mu-\nu+\rho)_k}{(1-\nu+\rho)_k k!} z^k \\
        &\hspace{10mm}-(-1)^{\nu-\rho} \frac{\Gamma(\mu+\rho)}{\Gamma(\rho) \Gamma(\nu-\rho+1)\Gamma(\mu-\nu+\rho)} \pFq{2}{1}{\mu,\nu}{\nu-\rho+1}{z} \,\log(z) \\
        &\hspace{10mm}+ (-1)^{\nu-\rho}\frac{\Gamma(\mu+\rho)}{\Gamma(\rho) \Gamma(\nu-\rho+1) \Gamma(\mu-\nu+\rho)} \sum_{k=0}^{\infty} \frac{(\mu)_k (\nu)_k}{(\nu-\rho+1)_k k!}\\
        &\hspace{24mm}\times\left[\psi(k+1)-\psi(k+\mu)+\psi(k+\nu+1-\rho)-\psi(k+\nu)\right] z^k ,
    \end{split}
\end{equation}
for $|z|<1$, $|\mathrm{arg}(z)|<\pi$, $(\nu-\rho)\in\mathbb{Z}^+$,
$\mathrm{Re}(\rho)>0$, $\rho\notin\mathbb{Z}^+$ and $\mathrm{Re}(\rho-\nu+\mu)>0$. 

Letting $\rho=\nu-n$ for $n\in\mathbb{Z}^+$ and performing the shift in variable from $z$ to $1-z$, the transformation equation \eqref{repcase4bx} assumes the form
\begin{equation}
\begin{split}
&\pFq{2}{1}{\mu,\nu}{\mu+\nu-n}{z}= \frac{\Gamma(\mu+\nu-n) (n-1)!}{\Gamma(\nu) \Gamma(\mu)} (1-z)^{-n} \sum_{k=0}^{n-1} \frac{(\nu-n)_k (\mu-n)_k}{(1-n)_k k!}(1-z)^k \\
&\hspace{14mm} + (-1)^n \frac{\Gamma(\mu+\nu-n)}{\Gamma(\nu-n)\Gamma(\mu-n) n!} \sum_{k=0}^{\infty} \frac{(\mu)_k (\nu)_k}{(n+1)_k k!} (1-z)^k \\
&\hspace{16mm} \times \left[-\log(1-z)+\psi(k+1)+\psi(k+n+1)-\psi(k+\mu)-\psi(k+\nu)\right]  ,
\end{split}
\end{equation}
which is valid for $|1-z|<1$, $|\mathrm{arg}(1-z)|<\pi$, $\nu\neq 0$, $\mathrm{Re}(\nu)>n$ and $\mathrm{Re}(\mu)>n$. For a fixed given $n$, the restrictions for $\nu$ and $\mu$ may be lifted by analytic continuation. The result is a known transformation equation \cite{wolfram}.

\subsubsection{Case $(\rho-\nu)\in\mathbb{Z}^+_0$} For this case all integrals are finite part integrals and the expansion assumes the form
\begin{equation}
\begin{split}
 &\int_0^{\infty} \frac{x^{\nu-1} (a+x)^{-\mu}}{(b+x)^{\rho}} \mathrm{d}x =  \sum_{k=0}^{\infty} {- \rho \choose k} \, b^{k} \, \bbint{0}{\infty} \frac{(a+x)^{-\mu}}{x^{k+\rho-\nu+1}} \,\mathrm{d}x
\\
& \hspace{24mm}- (-1)^{\nu-\rho}\frac{\sin(\pi\rho)}{\pi} \bbint{0}{b} \frac{x^{\nu-1} \, (a-x)^{-\mu}}{(b-x)^{\rho}} \, \ln(x)\, \mathrm{d}x,
\label{case4c}
\end{split}
\end{equation}
under the same conditions as for \eqref{case4b}. First, let us consider the case $\rho\neq\nu$. Substituting the finite part integrals \eqref{FPIpole} in the first term and performing hypergeometric summation, the first term becomes
\begin{equation}
\begin{split}
&\sum_{k=0}^{\infty} {- \rho \choose k} \, b^{k} \, \bbint{0}{\infty} \frac{(a+x)^{-\mu}}{x^{k+\rho-\nu+1}} \,\mathrm{d}x \\
&\hspace{12mm}= \frac{(-1)^{\rho-\nu}}{a^{\rho-\nu+\mu}} \frac{\Gamma(\mu+\rho-\nu)}{\Gamma(\mu) \Gamma(\rho-\nu+1)} \pFq{2}{1}{\rho,\mu+\rho-\nu}{\rho-\nu+1}{\frac{b}{a}} \,\ln(a) \\
&\hspace{12mm} + \frac{(-1)^{\rho-\nu}}{a^{\rho-\nu+\mu}} \frac{\Gamma(\mu+\rho-\nu)}{\Gamma(\mu) \Gamma(\rho-\nu+1)} \sum_{k=0}^{\infty} \frac{(\rho)_k (\mu+\rho-\nu)_k}{(\rho-\nu+1)_k k!} \\
&\hspace{28mm}\times\left[\psi(k+\rho-\nu+1)-\psi(k+\rho-\nu+\mu)\right] \left(\frac{b}{a}\right)^k,
\label{case4cx}
\end{split}
\end{equation}
which is valid for $a>b>0$, $(\rho-\nu)\in\mathbb{Z}^+$ and $\mathrm{Re}(\rho-\nu+\mu)>0$.

Since $(\rho-\nu)$ is a positive integer, the argument $(1+k+\nu-\rho)$  may be a positive or a negative integer depending on the value of the index $k$ in equations \eqref{fpingteg1x} and \eqref{fpingteg1xxx}. Then for $k=0, \dots,\rho-\nu-1$, $1/\Gamma(1+k+\nu-\rho)=0$, we have the vanishing finite parts 
\begin{equation}\label{fpingteg1x2}
\bbint{0}{1}\frac{y^{\nu+k-1}}{(1-y)^{\rho}}\,\mathrm{d}y=0,
\end{equation}
for $k=0, \dots, (\rho-\nu-1)$, $\mathrm{Re}(\nu)>0$, $\mathrm{Re}(\rho)>0$ and $\rho\notin\mathbb{Z}^+$. On the other hand, for $k=(\rho-\nu), (\rho-\nu)+1, \dots$, the argument of the gamma function is positive, so that it yields the same value given by equation \eqref{fpingteg1x},
\begin{equation}\label{be2}
\bbint{0}{1}\frac{y^{k+\nu-1}}{(1-y)^{\rho}}\,\mathrm{d} y = \frac{\pi \Gamma(k+\nu)}{\sin(\pi\rho) \Gamma(\rho) \Gamma(k+\nu-\rho+1)},
\end{equation} 
also for $\mathrm{Re}(\nu)>0$, $\mathrm{Re}(\rho)>0$ and $\rho\notin\mathbb{Z}^+$. The summation in equation \eqref{fpixxx} then starts at $k=\rho-\nu$, yielding the finite part integral
\begin{equation}\label{xfpi1}
\bbint{0}{1} \frac{y^{\nu-1} (1-\frac{b}{a}y)^{-\mu}}{(1-y)^{\rho}}\,\mathrm{d}y = \frac{\pi \Gamma(\rho-\nu+\mu)}{\sin(\pi \rho) \Gamma(\mu) \Gamma(\rho-\nu+1)} \left(\frac{b}{a}\right)^{\rho-\nu} \pFq{2}{1}{\mu,\rho}{\rho-\nu+1}{\frac{b}{a}} ,
\end{equation}
for $a>b>0$, $(\rho-\nu)\in\mathbb{Z}^+$,  $\mathrm{Re}(\nu)>0$, $\mathrm{Re}(\rho)>0$ and $\rho\notin\mathbb{Z}^+$, where a hypergeometric summation has been performed on the resulting infinite series. 

For the logarithmic integral, the value of the finite part depends on the value of $k$ as well. From equation \eqref{fpingteg1xx3}, we have the value 
\begin{equation}\label{fpingteg1xx2}
\begin{split}
\bbint{0}{1}\frac{y^{\nu+k-1}}{(1-y)^{\rho}}\,\ln(y)\,\mathrm{d}y = (-1)^{k+\rho-\nu+1} \frac{\pi\Gamma(\nu+k) (\rho-\nu-k-1)!}{\sin(\pi\rho) \Gamma(\rho)},
\end{split}
\end{equation}
for $k=0, \dots,(\rho-\nu-1)$; while the same value given by equation \eqref{fpingteg1xxx2} holds for $k=(\rho-\nu), \,(\rho-\nu)+1, \dots$,
\begin{equation}\label{fpingteg1xx22}
\bbint{0}{1}\frac{y^{\nu+k-1}}{(1-y)^{\rho}}\,\ln(y)\,\mathrm{d}y = \frac{\pi\Gamma(\nu+k)}{\sin(\pi\rho) \Gamma(\rho) \Gamma(1+k+\nu-\rho)} \left(\psi(k+\nu)-\psi(1+k+\nu-\rho)\right).
\end{equation} Both equations \eqref{fpingteg1xx2} and \eqref{fpingteg1xx22} are valid for $\mathrm{Re}(\nu)>0$, $\mathrm{Re}(\rho)>0$ and $\rho\notin\mathbb{Z}^+$. Then the finite part integral assumes the form
\begin{equation}\label{xfpi2}
\begin{split}
&\bbint{0}{1} \frac{y^{\nu-1} (1-\frac{b}{a}y)^{-\mu}}{(1-y)^{\rho}}\,\ln(y)\,\mathrm{d}y\\ &\hspace{14mm}=(-1)^{\rho-\nu} \frac{\pi \Gamma(\nu)}{\sin(\pi\rho) \Gamma(\rho)} \sum_{k=0}^{\rho-\nu-1} \frac{(\mu)_k (\nu)_k (\rho-\nu-k-1)!}{k!}\left(-\frac{b}{a}\right)^k\\
&\hspace{18mm} +\frac{\pi \Gamma(\rho-\nu+\mu)}{\sin(\pi\rho) \Gamma(\mu) \Gamma(\rho-\nu+1)} \left(\frac{b}{a}\right)^{\rho-\nu} \sum_{k=0}^{\infty} \frac{(\rho-\nu+\mu)_k (\rho)_k}{(\rho-\nu+1)_k k!} \\
&\hspace{34mm}\times \left[\psi(k+\rho)-\psi(k+1)\right]\left(\frac{b}{a}\right)^k
\end{split}
\end{equation}
which is valid for all $a>b>0$, $(\rho-\nu)\in\mathbb{Z}^+$, $\mathrm{Re}(\nu)>0$, $\mathrm{Re}(\rho)>0$ and $\rho\notin\mathbb{Z}^+$.

Substituting back the finite part integrals \eqref{xfpi1} and \eqref{xfpi2} into equation \eqref{integralorigx} yields the finite part integral
\begin{equation}\label{caca}
\begin{split}
&\frac{a^{\mu}}{b^{\nu-\rho}} \bbint{0}{b} \frac{x^{\nu-1} \, (a-x)^{-\mu}}{(b-x)^{\rho}} \, \ln(x)\, \mathrm{d}x\\
&\hspace{20mm}=\frac{\pi \Gamma(\rho-\nu+\mu)}{\sin(\pi \rho) \Gamma(\mu) \Gamma(\rho-\nu+1)} \left(\frac{b}{a}\right)^{\rho-\nu} \pFq{2}{1}{\mu,\rho}{\rho-\nu+1}{\frac{b}{a}}\,\ln(b)\\
&\hspace{22mm}+(-1)^{\rho-\nu} \frac{\pi \Gamma(\nu)}{\sin(\pi\rho) \Gamma(\rho)} \sum_{k=0}^{\rho-\nu-1} \frac{(\mu)_k (\nu)_k (\rho-\nu-k-1)!}{k!}\left(-\frac{b}{a}\right)^k\\
&\hspace{22mm} +\frac{\pi \Gamma(\rho-\nu+\mu)}{\sin(\pi\rho) \Gamma(\mu) \Gamma(\rho-\nu+1)} \left(\frac{b}{a}\right)^{\rho-\nu} \sum_{k=0}^{\infty} \frac{(\rho-\nu+\mu)_k (\rho)_k}{(\rho-\nu+1)_k k!} \\
&\hspace{34mm}\times \left[\psi(k+\rho)-\psi(k+1)\right]\left(\frac{b}{a}\right)^k,
\end{split}
\end{equation} 
for $a>b>0$, $(\rho-\nu)\in\mathbb{Z}^+$, $\mathrm{Re}(\nu)>0$, $\mathrm{Re}(\rho)>0$ and $\rho\notin\mathbb{Z}^+$. Gathering all the finite part integrals, the logarithmic terms again combine and the resulting expression is also in the variable $b/a$. Replacing $b/a$ with complex $z$, we obtain
\begin{equation}\label{resultx}
\begin{split}
& \pFq{2}{1}{\mu,\nu}{\mu+\rho}{1-z}
=(-1)^{\nu-\rho+1} \frac{\Gamma(\mu+\rho)}{\Gamma(\mu)\Gamma(\nu)\Gamma(\rho-\nu+1)} z^{\rho-\nu} \sum_{k=0}^{\infty} \frac{(\rho)_k (\mu+\rho-\nu)_k}{(\rho-\nu+1)_k k!} \\
&\hspace{12mm} \times \left[\log(z)+\psi(k+\rho-\nu+1)-\psi(k+\rho-\nu+\mu) -\psi(k+\rho)+\psi(k+1))\right] z^k \\
&\hspace{18mm} +\frac{\Gamma(\mu+\rho)}{\Gamma(\rho) \Gamma(\mu+\rho-\nu)} \sum_{k=0}^{\rho-\nu-1} (-1)^k \frac{(\mu)_k (\nu)_k(\rho-\nu-k-1)!}{k! } z^k ,
\end{split}
\end{equation}
for $|z|<1$ with $|\mathrm{arg}(z)|<\pi$, $(\rho-\nu)\in\mathbb{Z}^+$, 
$\mathrm{Re}(\nu)>0$, $\mathrm{Re}(\rho)>0$, $\rho\notin\mathbb{Z}^+$ and $\mathrm{Re}(\rho-\nu+\mu)>0$. 

Known representations of the Gauss hypergeometric function can be extracted from equation \eqref{resultx}.  Also shifting variable from $z$ to $1-z$ and writing $\rho=\nu+m$ for $m\in\mathbb{Z}^+$, we obtain
\begin{equation}\label{repcase4d}
\begin{split}
&\pFq{2}{1}{\mu,\nu}{\mu+\nu+m}{z} = \frac{\Gamma(\mu+\nu+m)}{\Gamma(\mu+m) \Gamma(\nu+m)} \sum_{k=0}^{m-1} \frac{(\mu)_k (\nu)_k (m-k-1)!}{k!} (z-1)^k \\
&\hspace{8mm}-(z-1)^m \frac{\Gamma(\mu+\nu+m)}{\Gamma(\nu)\Gamma(\mu) m!} \sum_{k=0}^{\infty} \frac{(\nu+m)_k (\mu+m)_k}{(m+1)_k k!} (1-z)^k \\ &\hspace{12mm}\times\left[\log(1-z)-\psi(k+m+1)+\psi(k+\mu+m)+\psi(k+\nu+m)-\psi(k+1)\right],
\end{split}
\end{equation}
for $|1-z|<1$, $|\mathrm{arg}(1-z)|<\pi$, $m\in\mathbb{Z}^+$, $\mathrm{Re}(\nu)>-m$ and $\mathrm{Re}(\mu)>-m$. For a given fixed $m$, the restrictions on $\nu$ and $\mu$ may be lifted by analytic continuation. The result reproduces the transformation tabulated in \cite[Eqn. 15.8.10]{dlmf} except that the transformation is for the corresponding regularized hypergeometric function there. 

Finally, we consider the case $\rho=\nu$. Under this condition, the Stieltjes integral evaluate to \eqref{case4c} with $\nu=\rho$, with all the finite part integrals involved as special values of \eqref{firstterm} and \eqref{coco} for $\nu=\rho$. Then the following transformation equation can be readily established,
\begin{equation}\label{repcase4c}
\begin{split}
&   \pFq{2}{1}{\mu,\rho}{\mu+\rho}{z}
=- \frac{\Gamma(\mu+\rho)}{\Gamma(\rho)\Gamma(\mu)}\, \pFq{2}{1}{\mu,\rho}{1}{z} \,\log(1-z)\\
&\hspace{10mm}+ \frac{\Gamma(\mu+\rho)}{\Gamma(\rho)\Gamma(\mu)}\, \sum_{k=0}^{\infty} \frac{(\mu)_k (\nu)_k}{(k!)^2}
\left[2\psi(k+1)-\psi(k+\mu)-\psi(k+\rho)\right] (1-z)^k , 
\end{split}
\end{equation}
valid for all $|1-z|<1$ and $|\mathrm{arg}(1-z)|<\pi$, $\mathrm{Re}(\mu)>0$, $\mathrm{Re}(\rho)>0$ and $\rho\notin\mathbb{Z}^+$. This can be taken as a special value of equation \eqref{resultx} under the same conditions provided we interpret an empty sum as zero. 

\section{Transformation equations arising for the Generalized Hypergeometric Function $_3F_2$}\label{section3f2}
We now wish to apply finite part integration on the integral representation \eqref{general} for the specific case of the generalized hypergeometric function $_3F_2$ for the family of parameters $\alpha_2=\{\nu,1\}$ and $\rho_1=\{n\}$, where $\nu\notin\mathbb{Z}$ and $n\in\mathbb{Z}^+$. Then from \eqref{general}, we obtain the representation 
\begin{equation}\label{3f2}
\pFq{3}{2}{\beta,\nu,1}{\beta+\sigma,n}{z} = \frac{\Gamma(\beta+\alpha)}{\Gamma(\beta) \Gamma(\sigma)} \int_0^{\infty} \frac{s^{\beta-1}}{(s+1)^{\beta+\alpha}}\, \pFq{2}{1}{\nu,1}{n}{\frac{z s}{s+1}}\,\mathrm{d}s ,
\end{equation}
for $\mathrm{Re}(\beta)>0$, $\mathrm{Re}(\sigma)>0$ and $|z|<1$, where
\begin{equation}\label{keykey}
\pFq{2}{1}{\nu,1}{n}{z} = \frac{(\nu-n)(n-1)!}{(\nu-n)_n z^{n-1}} \left[ (1-z)^{n-\nu-1} - \sum_{k=0}^{n-2} \frac{(\nu-n+1)_{k}}{k!} z^k\right] .
\end{equation}
The identity \eqref{keykey} is established in the Appendix. Substituting \eqref{keykey} back into \eqref{3f2} and distributing the integration yield
\begin{equation}\label{3f2x}
\begin{split}
&\frac{\Gamma(\beta) \Gamma(\sigma) (\nu-n)_n\, z^{n-1}}{\Gamma(\beta+\alpha) (\nu-n)(n-1)!}\,\pFq{3}{2}{\beta,\nu,1}{\beta+\sigma,n}{z} =  \int_0^{\infty} \frac{s^{\beta-n} (1+(1-z) s)^{n-\nu-1}}{(s+1)^{\beta+\alpha-\nu}}\,\mathrm{d}s \\
&\hspace{26mm} -\sum_{k=0}^{n-2} \frac{(\nu-n+1)_{k}}{k!} z^k \int_0^{\infty} \frac{s^{\beta+k-n}}{(1+s)^{\beta+\sigma+k-n+1}} \,\mathrm{d}s .
\end{split} 
\end{equation}
In order for the distribution of the integration to be valid, the condition on $\beta$ must now be modified to $\mathrm{Re}(\beta)>(n-1)$, with all other conditions unchanged. We will later relax this condition by analytic continuation. The second integral is tabulated and is given by
\begin{equation}\label{leadint}
\int_0^{\infty} \frac{s^{\beta+k-n}}{(1+s)^{\beta+\sigma+k-n+1}} \,\mathrm{d}s = \frac{\Gamma(k-n+\beta+1)}{\Gamma(k-n+\beta+\sigma+1)}
\end{equation}
for $k+\mathrm{Re}(\beta)>(n-1)$ and $\mathrm{Re}(\sigma)>0$. Under the condition $\mathrm{\beta}>(n-1)$, equation \eqref{leadint} holds for all $k$. 

We will derive consequences of the representation \eqref{3f2x} by performing finite part integration on the integral
\begin{equation}\label{integ3f2}
\int_0^{\infty} \frac{s^{\beta-n} (1+(1-z) s)^{n-\nu-1}}{(s+1)^{\beta+\sigma-\nu}}\,\mathrm{d}s,
\end{equation}
for $n\in\mathbb{Z}^+$, $\mathrm{Re}(\beta)>(n-1)$, $\nu\notin\mathbb{Z}$ and $\mathrm{Re}(\sigma)>0$. We interpret \eqref{integ3f2} as a Stieltjes transform of the function $s^{\beta-n} (1+(1-z) s)^{n-\nu-1}$ with the kernel of transformation given by $(s+1)^{-\beta-\sigma+\nu}$, and apply the scheme developed above to evaluate the integral. With the choice of kernel of transformation, we impose the condition $\mathrm{Re}(\beta+\sigma-\nu)>0$. The fundamental divergent integrals are identified by expanding the kernel,
\begin{equation}\label{ex}
\frac{1}{(s+1)^{\beta+\alpha-\nu}} = \sum_{k=0}^{\infty} \binom{-\beta-\sigma+\nu}{k} \frac{1}{s^{\beta+\sigma-\nu+k}}
\end{equation}
inside the integral and distributing the integration. The divergent integrals are
\begin{equation}\label{fun}
\int_0^{\infty} \frac{(1+(1-z) s)^{n-\nu-1}}{s^{n+\sigma-\nu+k}}\,\mathrm{d}s 
\end{equation}
for sufficiently large values of $k$. Under the condition $\mathrm{Re}(\sigma)>0$,  integral \eqref{fun} is integrable at infinity and the divergence arises at the origin only. Notice that the nature of the singularity at the origin is independent of $\beta$. Since $n$ and $k$ are both integers, whether the singularity at the origin is a pole or branch point depends on the difference $\sigma-\nu$. 

As for the Gauss function, there are two distinct cases: $(\sigma-\nu)\notin\mathbb{Z}$ for a branch point singularity at the origin, and $(\sigma-\nu)\in\mathbb{Z}$ for a pole singularity. The later further splits into the case of $(n+\sigma-\nu)$ a natural number and  $(n+\sigma-\nu)$ a negative integer. A couple of other cases arise depending on whether $(\beta+\sigma-\nu)$ is an integer or not. The integer case leads to a pole singularity in the kernel at the point $z=-1$. We will only consider the non-integer case. Again to facilitate the process of finite part integration, we assume that $0<z<1$ and let $z=x$, with the intention to eventually appeal to analytic continuation to lift the restriction. To make use of the results obtained for the Gauss function, we recast integral \eqref{integ3f2} by factoring out $(1-x)$ in \eqref{integ3f2} and perform finite integration on the integral
\begin{equation}\label{integ3f2x}
\int_0^{\infty} \frac{s^{\beta-n} ((1-x)^{-1}+ s)^{n-\nu-1}}{(s+1)^{\beta+\sigma-\nu}}\,\mathrm{d}s ,
\end{equation}
under the same conditions for \eqref{integ3f2}.

The complex extension of $s^{\beta-n} ((1-x)^{-1}+ s)^{n-\nu-1}$, denoted by $z^{\beta-n} ((1-x)^{-1}+ z)^{n-\nu-1}$, has either a pole or a branch point at $z=-(1-x)^{-1}$, depending on whether $(n-\nu-1)$ is a negative integer (a pole) or a non-integer (a branch point). When $-(1-x)^{-1}$ happens to be a branch point, we choose the branch cut to be $(-\infty,-(1-x)^{-1}]$, in accordance with the general requirements in Section-\ref{finitepartintegration}.  But again the final result is independent of the nature of the singularity at $-(1-x)^{-1}$ so that no restriction is imposed on $\nu$ at this point. Under the condition $\mathrm{Re}(\beta)>(n-1)$, the origin can only be a zero or a branch point depending on $\beta$. When $\beta$ is a non-integer, we choose the branch cut to be $[0,\infty)$. On the other hand, the complex extension of the kernel  $(s+1)^{-\beta-\sigma+\nu}$, denoted by $(z+1)^{-\beta-\sigma+\nu}$, has a pole or branch point singularity at $-1$. If $-1$ is a branch point, the branch cut is again to chosen to be $[-1,\infty)$. Since $0<x<1$, the singularity of $((1-x)^{-1}+ z)^{n-\nu-1}$ is to the right of the singularity of the kernel $(z+1)^{-\beta-\sigma+\nu}$. This allows our scheme of finite part integration to be applied in the evaluation of the integral \eqref{integ3f2x}. 

\subsection{Case: $(\sigma-\nu)\notin\mathbb{Z}$ and $(\beta+\sigma-\nu)\notin\mathbb{Z}$}
Under this case, the origin is a branch-point singularity of the divergent integral \eqref{fun} and the kernel has a branch point singularity at $-1$. First we assume that $(\sigma-\nu)\geq 0$ so that the integral \eqref{fun} is divergent for all $k=0, 1, 2, \dots$.   Then we extract the integral \eqref{integ3f2x} from the contour integral
\begin{equation}
\frac{1}{e^{-2\pi i (\sigma-\nu)}-1} \int_{\mathrm{C}} \frac{z^{\beta-n} ((1-x)^{-1}+z)^{n-\nu-1}}{(1+z)^{\beta+\sigma-\nu}}\,\mathrm{d}z . 
\end{equation}
Again collapsing the contour $\mathrm{C}$ into the contour $\mathrm{C}'$, the integral emmerges to take the representation
\begin{equation}\label{integ3f2xx}
\begin{split}
&\int_0^{\infty} \frac{s^{\beta-n} ((1-x)^{-1}+ s)^{n-\nu-1}}{(s+1)^{\beta+\alpha-\nu}}\,\mathrm{d}s = \frac{1}{e^{-2\pi i (\sigma-\nu)}-1} \int_{\mathrm{C}} \frac{z^{\beta-n} ((1-x)^{-1}+z)^{n-\nu-1}}{(1+z)^{\beta+\sigma-\nu}}\,\mathrm{d}z\\
&\hspace{18mm} - \frac{\sin(\pi(\beta+\sigma-\nu))}{\sin(\pi(\sigma-\nu))} \lim_{\epsilon\rightarrow 0} \left[\int_0^{1-\epsilon} \frac{s^{\beta-n} ((1-x)^{-1} -s)^{n-\nu-1}}{(1-s)^{\beta+\sigma-\nu}}\,\mathrm{d}s\right. \\
&\hspace{18mm}\left. +\frac{\sin(\pi (\sigma-\nu))}{\sin(\pi(\beta+\sigma-\nu)) (e^{-2\pi i (\sigma-\nu)}-1)}\int_{\mathrm{C}_{\epsilon}} \frac{z^{\beta-n} ((1-x)^{-1}+z)^{n-\nu-1}}{(1+z)^{\beta+\sigma-\nu}}\,\mathrm{d}z \right].
\end{split}
\end{equation}
The limit can again be established as a finite part integral. 
Expanding the first term and distributing the integration, the resulting contour integrals are again identified to be the finite parts of the divergent integrals \eqref{fun}. Then the integral evaluates to
\begin{equation}\label{keyx}
\begin{split}
&\int_0^{\infty} \frac{s^{\beta-n} ((1-x)^{-1}+ s)^{n-\nu-1}}{(s+1)^{\beta+\sigma-\nu}}\,\mathrm{d}s = \sum_{k=0}^{\infty} \binom{-\beta-\sigma+\nu}{k} \bbint{0}{\infty} \frac{((1-x)^{-1} + s)^{n-\nu-1}}{s^{n+\sigma-\nu+k}}\,\mathrm{d}s \\
& \hspace{34mm}- \frac{\sin(\pi(\beta+\sigma-\nu))}{\sin(\pi(\sigma-\nu))}  \bbint{0}{1} \frac{s^{\beta-n} ((1-x)^{-1} -s)^{n-\nu-1}}{(1-s)^{\beta+\sigma-\nu}}\,\mathrm{d}s,
\end{split}
\end{equation}
which is valid for $0<x<1$, $n\in\mathbb{Z}^+$, $(\sigma-\nu)\notin\mathbb{Z}$, $\nu\notin\mathbb{Z}$, $\mathrm{Re}(\beta)>(n-1)$, $(\beta+\sigma-\nu)\notin\mathbb{Z}$, $\mathrm{Re}(\beta+\sigma-\nu)>0$ and $\mathrm{Re}(\sigma)>0$. 

The fundamental finite part integrals in equation \eqref{keyx} follows from equation \eqref{finxx} with the substitutions $\lambda\rightarrow (n +\sigma-\nu+k)$, $\upsilon\rightarrow(\nu+1-n)$, and $s\rightarrow (1-x)^{-1}$. The desired finite part integral is given by
\begin{equation}
\bbint{0}{\infty} \frac{((1-x)^{-1} + s)^{n-\nu-1}}{s^{n+\sigma-\nu+k}}\,\mathrm{d}s =  \frac{\pi (-1)^{n+k}\Gamma(\sigma) (\sigma)_k(1-x)^{\sigma+k} }{\sin(\pi(\sigma-\nu)) \Gamma(1+\nu-n) \Gamma(\sigma-\nu+n) (\sigma-\nu+n)_k},
\end{equation}
for all $k\in\mathbb{Z}_0^+$ and $\mathrm{Re}(\sigma)>0$. Substitution of this value yields for the first term of \eqref{keyx},  
\begin{equation}\label{ft3f2}
\begin{split}
&\sum_{k=0}^{\infty} \binom{-\beta-\sigma+\nu}{k} \bbint{0}{\infty} \frac{((1-x)^{-1} + s)^{n-\nu-1}}{s^{n+\sigma-\nu+k}}\,\mathrm{d}s\\
& \hspace{12mm}= \frac{\pi (-1)^n \Gamma(\sigma) (1-x)^{\sigma}}{\sin(\pi(\sigma-\nu)) \Gamma(1+\nu-n) \Gamma(\sigma-\nu+n)}  \pFq{2}{1}{\sigma,\beta+\sigma-\nu}{\sigma-\nu+n}{1-x} ,
\end{split}
\end{equation}
for $1>x>0$, $(\sigma-\nu)\notin\mathbb{Z}$, $\nu\notin\mathbb{Z}$ and $\mathrm{Re}(\sigma)>0$, where a hypergeometric summation has been performed to arrive at the right hand side. 

For the progenic finite part integral, we use the same method leading to equations \eqref{fpixxx} and \eqref{fpixxx2}. We factor $(1-x)^{-1}$ out from the integral to cast the finite-part integral into the form of equation \eqref{fpixxx} along with the identifications $\nu\rightarrow (\beta-n+1)$, $b/a\rightarrow (1-x)$ and $\rho\rightarrow(\beta+\sigma-\nu)$. We have
\begin{equation}\label{boo}
\begin{split}
&\bbint{0}{1} \frac{s^{\beta-n} ((1-x)^{-1}-s)^{n-\nu-1}}{(1-s)^{\beta+\sigma-\nu}}\,\mathrm{d}s\\
&\hspace{12mm} = (1-x)^{-(n-\nu-1)} \sum_{k=0}^{\infty} (-1)^k \binom{n-\nu-1}{k} (1-x)^k \bbint{0}{1} \frac{s^{(k+\beta-n+1)-1}}{(1-s)^{\beta+\sigma-\nu}}\,\mathrm{d}s . 
\end{split}
\end{equation}
for $1>x>0$, $\mathrm{Re}(\beta)>(n-1)$, $\mathrm{Re}(\beta+\sigma-\nu)>0$ and $(\beta+\sigma-\nu)\notin\mathbb{Z}$. The finite part integral in the right hand side is just \eqref{fpingteg1x} with $\sigma\rightarrow k+\beta-n+1$ and $\rho\rightarrow\beta+\sigma-\nu$. This yields the finite part integral
\begin{equation}\label{boox}
\bbint{0}{1} \frac{s^{(k+\beta-n+1)-1}}{(1-s)^{\beta+\sigma-\nu}}\,\mathrm{d}s = \frac{\pi \Gamma(\beta-n+1) (\beta-n+1)_k}{\sin(\pi(\beta+\sigma-\nu)) \Gamma(\beta+\sigma-\nu) \Gamma(2-n-\sigma+\nu) (2-n-\sigma+\nu)_k}
\end{equation}
for all $k\in\mathbb{Z}^+_0$, $\mathrm{Re}(\beta)>(n-1)$, $\mathrm{Re}(\beta+\sigma-\nu)>0$ and $(\beta+\sigma-\nu)\notin\mathbb{Z}^+$. Substituting \eqref{boox} back into \eqref{boo} and performing hypergeometric summation yield
\begin{equation}\label{st3f2}
\begin{split}
&\bbint{0}{1} \frac{s^{\beta-n} ((1-x)^{-1}-s)^{n-\nu-1}}{(1-s)^{\beta+\sigma-\nu}}\,\mathrm{d}s = \frac{\pi \Gamma(\beta-n+1) (1-x)^{-(n-\nu-1)}}{\sin(\pi(\beta+\sigma-\nu)) \Gamma(\beta+\sigma-\nu) \Gamma(2-n-\sigma+\nu)}\\
&\hspace{54mm}\times \pFq{2}{1}{1+\nu-n,\beta-n+1}{2-n-\sigma+\nu}{1-x},
\end{split}
\end{equation}
under the conditions for \eqref{boo}. The desired finite part integral \eqref{keyx} can now obtained by substituting \eqref{ft3f2} and \eqref{st3f2} back into \eqref{keyx}. The result is
\begin{equation}\label{keyxx}
\begin{split}
&\int_0^{\infty} \frac{s^{\beta-n} ((1-x)^{-1}+ s)^{n-\nu-1}}{(s+1)^{\beta+\sigma-\nu}}\,\mathrm{d}s\\
& \hspace{12mm}=  \frac{\pi (-1)^n  \Gamma(\sigma) (1-x)^{\sigma}}{\sin(\pi(\sigma-\nu)) \Gamma(1+\nu-n) \Gamma(\sigma-\nu+n)}  \pFq{2}{1}{\sigma,\beta+\sigma-\nu}{\sigma-\nu+n}{1-x}\\
& \hspace{12mm} - \frac{\pi \Gamma(\beta-n+1) (1-x)^{-(n-\nu-1)}}{\sin(\pi(\sigma-\nu)) \Gamma(\beta+\sigma-\nu) \Gamma(2-n-\sigma+\nu)} \pFq{2}{1}{1+\nu-n,\beta-n+1}{2-n-\sigma+\nu}{1-x} ,
\end{split}
\end{equation}
valid for $0<x<1$, $n\in\mathbb{Z}^+$, $\mathrm{Re}(\sigma)>0$, $\mathrm{Re}(\beta)>(n-1)$, $\nu\notin\mathbb{Z}$ and $(\sigma-\nu)\notin\mathbb{Z}$.

Substituting equations \eqref{ft3f2} and \eqref{keyx} back into equation \eqref{3f2x} we obtain the following identity,
\begin{equation}\label{iden}
\begin{split}
&\frac{\Gamma(\beta) (\nu-n)_n  }{\Gamma(\beta+\sigma)(\nu-n)(n-1)!}z^{n-1}\,\pFq{3}{2}{\beta,\nu,1}{\beta+\sigma,n}{x}=\frac{\pi}{\sin(\pi(\sigma-\nu))}\\
&\hspace{8mm} \times \left[\frac{(-1)^n (1-x)^{\sigma+n-\nu-1}}{\Gamma(1+\nu-n)\Gamma(\sigma-\nu+n)} \pFq{2}{1}{\sigma,\beta+\sigma-\nu}{\sigma-\nu+n}{1-x}\right.\\
&\hspace{8mm}-\left. \frac{\pi \Gamma(\beta-n+1)}{\Gamma(\sigma) \Gamma(\beta+\sigma-\nu) \Gamma(2-n-\sigma+\nu)} \pFq{2}{1}{1+\nu-n,\beta-\nu+1}{2-n-\sigma+\nu}{1-x}\right] \\
&\hspace{8mm} -  \sum_{k=0}^{n-2} \frac{(\nu-n+1)_{k} \Gamma(k-n+\beta+1)}{\Gamma(k-n+\beta+\sigma+1) k!} x^k , 
\end{split}
\end{equation}
under the same conditions as for \eqref{keyxx}. An empty sum in \eqref{iden} is interpreted as zero. This identity can be extended for complex $z$, with all functions involved given by their principal values. Implementing the extension in the complex plane and performing an expansion of the finite sum about $z=1$ yield the following transformation equation,
\begin{equation}
\begin{split}
&\frac{\Gamma(\beta) (\nu-n)_n  }{\Gamma(\beta+\sigma)(\nu-n)(n-1)!}z^{n-1}\, \pFq{3}{2}{\beta,\nu,1}{\beta+\sigma,n}{z}\\
&\hspace{6mm}=
\frac{(-1)^n \pi \csc(\pi(\sigma-\nu)) (1-z)^{\sigma+n-\nu-1}}{\Gamma(1+\nu-n)\Gamma(\sigma-\nu+n)} \pFq{2}{1}{\sigma,\beta+\sigma-\nu}{\sigma-\nu+n}{1-z}\\
&\hspace{8mm}-\frac{\pi \csc(\pi(\sigma-\nu))\Gamma(\beta-n+1)}{\Gamma(\sigma) \Gamma(\beta+\sigma-\nu) \Gamma(2-n-\sigma+\nu)} \pFq{2}{1}{1+\nu-n,\beta-\nu+1}{2-n-\sigma+\nu}{1-z} \\
&\hspace{8mm} - \sum_{l=0}^{n-2} \frac{1}{l!}\left[\sum_{k=l}^{n-2}\frac{(\nu-n+1)_{k}  \Gamma(k-n+\beta+1)k!}{\Gamma(k-n+\beta+\sigma+1) (k-l)!}\right] (z-1)^l
\end{split}
\end{equation}
for $|1-z|<1$,  $n\in\mathbb{Z}^+$, $(\sigma-\nu)\notin\mathbb{Z}$, $\nu\notin\mathbb{Z}$, $\mathrm{Re}(\beta)>(n-1)$, $(\beta+\sigma-\nu)\notin\mathbb{Z}$, $\mathrm{Re}(\beta+\sigma-\nu)>0$ and $\mathrm{Re}(\sigma)>0$. For a fixed given positive integer $n$, the restrictions on $\sigma$, $\nu$ and $\beta$ may be lifted by analytic continuation provided that $(\sigma-\nu)\notin\mathbb{Z}$.   The result appears to be new.

\subsection{Case $(\sigma-\nu)\in\mathbb{Z}$ and $\beta\notin\mathbb{Z}$} Under this condition, the divergence of \eqref{fun} at the origin is due to a pole singularity and the kernel has a branch point singularity. The integral \eqref{integ3f2x} is then extracted from the contour integral
\begin{equation}
\frac{1}{2\pi i} \int_{\mathrm{C}} \frac{z^{\beta-n} ((1-x)^{-1}+z)^{n-\nu-1}}{(1+s)^{\beta+\sigma-\nu}} (\log(z)-i\pi)\,\mathrm{d}z .
\end{equation}
On deforming the contour $\mathrm{C}$ to $\mathrm{C}'$, the integral \eqref{integ3f2x} emerges to assume the representation
\begin{equation}
\begin{split}\label{veve}
&\int_0^{\infty} \frac{s^{\beta-n} ((1-x)^{-1}+s)^{n-\nu-1}}{(1+s)^{\beta+\sigma-\nu}}\,\mathrm{d}s \\
&\hspace{16mm} = \frac{1}{2\pi i} \int_{\mathrm{C}} \frac{z^{\beta-n} ((1-x)^{-1}+z)^{n-\nu-1}}{(1+s)^{\beta+\sigma-\nu}} (\log(z)-i\pi)\,\mathrm{d}z\\
&\hspace{18mm}+(-1)^n\frac{\sin(\pi\beta)}{\pi} \lim_{\epsilon\rightarrow 0} \left[\int_0^{1-\epsilon} \frac{s^{\beta-n} ((1-x)^{-1} -s)^{n-\nu-1}}{(1-s)^{\beta+\sigma-\nu}} \ln(s)\,\mathrm{d}s\right. \\
&\hspace{18mm}\left. -\frac{(-1)^n}{2\pi i \sin(\pi\beta)}\int_{\mathrm{C}_{\epsilon}} \frac{z^{\beta-n} ((1-x)^{-1}+z)^{n-\nu-1}}{(1+z)^{\beta+\sigma-\nu}} (\log(z)-i\pi)\,\mathrm{d}z \right].
\end{split}
\end{equation} 
Again the second term is the finite part of the divergent integral $\int_0^1 s^{\beta-n} ((1-x)^{-1}-s)^{n-\nu-1} (1-s)^{-(\beta+\sigma-\nu)}\ln(s)\,\mathrm{d}s$. Expanding the kernel $(1+s)^{-(\beta+\sigma-\nu)}$ in the first term and performing term by term integration evaluate the integral in the form
\begin{equation}\label{fullint}
\begin{split}
&\int_0^{\infty} \frac{s^{\beta-n} ((1-x)^{-1}+s)^{n-\nu-1}}{(1+s)^{\beta+\sigma-\nu}}\,\mathrm{d}s = \sum_{k=0}^{\infty} \binom{-\beta-\sigma+\nu}{k} \bbint{0}{\infty} \frac{((1-x)^{-1}+s)^{n-\nu-1}}{s^{\sigma-\nu+n+k}}\,\mathrm{d}s \\
&\hspace{22mm} +(-1)^n\frac{\sin(\pi\beta)}{\pi} \bbint{0}{1} \frac{s^{\beta-n} ((1-x)^{-1} -s)^{n-\nu-1}}{(1-s)^{\beta+\sigma-\nu}} \ln(s)\,\mathrm{d}s,
\end{split}
\end{equation}
for $0<x<1$, $(\sigma-\nu)\in\mathbb{Z}$, $\nu\notin\mathbb{Z}$, $\mathrm{Re}(\beta)>(n-1)$, $\mathrm{Re}(\sigma)>0$, $\mathrm{Re}(\beta+\sigma-\nu)>0$ and $\beta\notin\mathbb{Z}$. There are two cases here: When $(\sigma-\nu+n)\in\mathbb{Z}^+$ and when $(\sigma-\nu+n)\in\mathbb{Z}^-_0$. In the former, all the integrals are finite part integrals; and in the later, the first few integrals are convergent integrals. We will take each separately. 

\subsubsection{Case $(\sigma-\nu+n)\in\mathbb{Z}^+$} The fundamental finite part integrals are given by equation \eqref{FPIpole} with the substitutions $\upsilon\rightarrow (\nu+1-n)$, $n\rightarrow (\sigma-\nu+n+k-1)$ and $s\rightarrow (1-x)^{-1}$. The result is
\begin{equation}\label{fpifun}
\begin{split}
\bbint{0}{\infty} \frac{((1-x)^{-1}+s)^{n-\nu-1}}{s^{\sigma-\nu+n+k}}\,\mathrm{d}s =& \frac{(-1)^{\sigma-\nu+n+k-1}(1+\nu-n)_{\sigma-\nu+n+k-1} (1-x)^{\sigma+k}}{(\sigma-\nu+n+k-1)!} \\
& \times \left[-\ln(1-x) + \psi(\sigma-\nu+n+k)-\psi(\sigma+k)\right],
\end{split}
\end{equation} 
for $1>x>0$ and $\mathrm{Re}(\sigma)>0$. Substituting \eqref{fpifun} back into the first term of \eqref{fullint}, we have
\begin{equation}\label{term1}
\begin{split}
 &\sum_{k=0}^{\infty} \binom{-\beta-\sigma+\nu}{k} \bbint{0}{\infty} \frac{((1-x)^{-1}+s)^{n-\nu-1}}{s^{\sigma-\nu+n+k}}\,\mathrm{d}s \\
 &\hspace{12mm}= \frac{(-1)^{\sigma-\nu+n-1} \Gamma(\sigma)}{\Gamma(1+\nu-n)}
 \sum_{k=0}^{\infty} \frac{(\beta+\sigma-\nu)_k (\sigma)_k}{(\sigma-\nu+n+k-1)! k!}\\
 &\hspace{20mm}\times \left[-\ln(1-x) + \psi(\sigma-\nu+n+k)-\psi(\sigma+k)\right] (1-x)^{\sigma+k}
\end{split}
\end{equation}
under the same conditions as for \eqref{fpifun}. 

Also for the progenic finite part integral, we perform an expansion in the same way we did to arrive at \eqref{fpixxx2}. We obtain the expansion
\begin{equation}\label{mofpi}
\begin{split}
	&\bbint{0}{1} \frac{s^{\beta-n} ((1-x)^{-1}-s)^{n-\nu-1}}{(1-s)^{\beta+\sigma-\nu}}\ln(s)\,\mathrm{d}s = (1-x)^{-(n-\nu-1)} \\
	&\hspace{34mm}\times \sum_{k=0}^{\infty} \frac{(1+\nu-n)_k}{k!} (1-x)^k \bbint{0}{1}\frac{s^{\beta+k-n} }{(1-s)^{\beta+\sigma-\nu}} \ln(s)\,\mathrm{d}s ,
	\end{split}
\end{equation}
for $0<x<1$, $\mathrm{Re}(\beta)>(n-1)$, $\mathrm{Re}(\beta+\sigma-\nu)>0$ and $\beta\notin\mathbb{Z}$. The finite part integral in the right hand side of \eqref{mofpi} is given by \eqref{fpingteg1xxx2} or \eqref{fpingteg1xx3}, with the substitutions $\sigma\rightarrow (\beta+k-n+1)$, $\rho\rightarrow (\beta+\sigma-\nu)$ and $(\sigma-\rho)\rightarrow (k-n+1-\sigma+\nu)$. Since $(\sigma-\nu+n)$ is a positive integer, the value of the finite part depends on whether $(k-n+1-\sigma+\nu)$ is a negative integer or not. We have the following values,
\begin{equation}\label{fpi1}
\begin{split}
&\bbint{0}{1}\frac{s^{\beta+k-n} }{(1-s)^{\beta+\sigma-\nu}} \ln(s)\,\mathrm{d}s= \frac{(-1)^{n-1-k}\pi}{\sin(\pi\beta)}\\
&\hspace{24mm} \times \frac{\Gamma(\beta-n+1) (\beta-n+1)_k (\sigma-\nu+n-2-k)!}{\Gamma(\beta+\sigma-\nu)}
\end{split}
\end{equation}
for $k=0,\dots, (\sigma-\nu+n-2)$, and
\begin{equation}\label{fpi2}
 \begin{split}
 &\bbint{0}{1}\frac{s^{\beta+k-n} }{(1-s)^{\beta+\sigma-\nu}} \ln(s)\,\mathrm{d}s = \frac{(-1)^{\sigma-\nu}\pi}{\sin(\pi\beta)}\\
 &\hspace{14mm} \times  \frac{\Gamma(\beta+k-n+1) (\psi(\beta+k-n+1)-\psi(k-n+2-\sigma+\nu))}{\Gamma(\beta+\sigma-\nu) \Gamma(k-n-\sigma+\nu+2)}
 \end{split}
\end{equation}
for $k=\sigma-\nu+n-1, \sigma-\nu+n, \dots$. Both equations \eqref{fpi1} and \eqref{fpi2} are valid for $\mathrm{Re}(\beta)>(n-1)$, $\mathrm{Re}(\beta+\sigma-\nu)>0$ and $\beta\notin\mathbb{Z}$. Substituting \eqref{fpi1} and \eqref{fpi2} back into \eqref{mofpi} yields the finite part integral
\begin{equation}\label{term2}
\begin{split}
	&\bbint{0}{1} \frac{s^{\beta-n} ((1-x)^{-1}-s)^{n-\nu-1}}{(1-s)^{\beta+\sigma-\nu}}\ln(s)\,\mathrm{d}s = \frac{\pi}{\sin(\pi\beta) \Gamma(\beta+\sigma-\nu)} (1-x)^{-(n-\nu-1)} \\
	&\hspace{10mm}\times \sum_{k=0}^{\sigma-\nu+n-2} \frac{(-1)^k \Gamma(k+\beta-n+1) (1+\nu-n)_k (\sigma-\nu+n-2-k)!}{k!} (1-x)^k \\
	&\hspace{10mm} + \frac{(-1)^{\sigma-\nu} \pi \Gamma(\sigma)}{\sin(\pi\beta) \Gamma(1+\nu-n)} \sum_{k=0}^{\infty} \frac{(\sigma)_k (\beta+\sigma-\nu)_k}{(k+\sigma-\nu+n-1)! k!}\\
	& \hspace{54mm}\times  \left[\psi(k+\beta+\sigma-\nu)-\psi(k+1)\right] (1-x)^{k+\sigma}
\end{split}
\end{equation}
which valid under the conditions for equation \eqref{mofpi}.  

The integral can now be evaluated. Substituting equations \eqref{term1} and \eqref{term2} back into equation \eqref{fullint}, we obtain
\begin{equation}\label{fullintx}
\begin{split}
&\int_0^{\infty} \frac{s^{\beta-n} ((1-x)^{-1}+s)^{n-\nu-1}}{(1+s)^{\beta+\sigma-\nu}}\,\mathrm{d}s =\frac{(-1)^n (1-x)^{-(n-\nu-1)}}{\Gamma(\beta+\sigma-\nu)}  \\
&\hspace{8mm}\times \sum_{k=0}^{\sigma-\nu+n-2} \frac{(-1)^k \Gamma(k+\beta-n+1) (1+\nu-n)_k (\sigma-\nu+n-2-k)!}{k!} (1-x)^k\\
&\hspace{8mm} + \frac{(-1)^{\sigma-\nu+n} \Gamma(\sigma)}{\Gamma(1+\nu-n)} \sum_{k=0}^{\infty} \frac{(\sigma)_k (\beta+\sigma-\nu)_k}{(k+\sigma-\nu+n-1)! k!} 
\left[\ln(1-x)-\psi(\sigma-\nu+n+k)\right. \\
&\hspace{18mm}\left. +\psi(k+\sigma)
 + \psi(k+\beta+\sigma-\nu)-\psi(k+1)\right] (1-x)^{k+\sigma},
\end{split}
\end{equation}
under the same conditions as in \eqref{mofpi}.

Finally, we gather equations \eqref{term1} and \eqref{fullintx} and substitute them back in the right hand side of \eqref{fpifun}.  The result, together with \eqref{leadint}, is returned to \eqref{3f2x}. We then extend the resulting expression by extending $x$ to complex values $z$ by analytic continuation. We  obtain
\begin{equation}\label{res2}
\begin{split}
&\frac{\Gamma(\beta) (\nu-n)_n  }{\Gamma(\beta+\sigma)(\nu-n)(n-1)!}z^{n-1}\, \pFq{3}{2}{\beta,\nu,1}{\beta+\sigma,n}{z}= \frac{1}{\Gamma(\sigma)\Gamma(\beta+\sigma-\nu)}\\
&\hspace{8mm}\times\sum_{k=0}^{\sigma-\nu+n-2}\frac{(-1)^k}{k!} \Gamma(k+\beta-n+1) (1+\nu-n)_k (\sigma-\nu+n-2-k)! (1-z)^k\\
&\hspace{8mm} + \frac{(-1)^{\sigma-\nu+n}}{\Gamma(\nu-n+1)} (1-z)^{\sigma+n-\nu-1} \sum_{k=0}^{\infty} \frac{(\sigma)_k (\beta+\sigma-\nu)_k}{(k+\sigma-\nu+n-1)! k!} (1-z)^k\\
&\hspace{8mm}\times \left[\log(1-z)-\psi(\sigma-\nu+n+k)+\psi(\sigma+k) + \psi(k+\beta+\sigma-\nu)-\psi(k+1)\right]\\
&\hspace{8mm} - \sum_{l=0}^{n-2} \left[\sum_{k=l}^{n-2}\frac{(\nu-n+1)_{k} \Gamma(k-n+\beta+1) k!}{\Gamma(k-n+\beta+\sigma+1) (k-l)!}\right] \frac{(z-1)^l}{l!}  
\end{split}
\end{equation}
for  all $|1-z|<1$, $n\in\mathbb{Z}^+$, $\nu\notin\mathbb{Z}$, $|\mathrm{arg}(1-z)|<\pi$, $(n+\sigma-\nu)\in\mathbb{Z}^+$, $\mathrm{Re}(\sigma)>0$,  $\mathrm{Re}(\beta+\sigma-\nu)>0$ and $\beta\notin\mathbb{Z}$, where we have performed the same expansion we did above to arrive at the last term of \eqref{res2}.  

We specialize for the case $n+\sigma-\nu=m\in\mathbb{Z}^+$ and eliminate $\sigma$ in favor of $n$, $m$, and $\nu$, i.e. $\sigma=m-n+\nu$. The result is
\begin{equation}\label{res2x}
\begin{split}
&\frac{\Gamma(\beta) (\nu-n)_n  }{\Gamma(\beta+m-n+\nu)(\nu-n)(n-1)!}z^{n-1}\, \pFq{3}{2}{\beta,\nu,1}{\beta+m-n+\nu,n}{z}\\
&\hspace{8mm}= \frac{1}{\Gamma(m-n+\nu)\Gamma(\beta+m-n)}
\sum_{k=0}^{m-2}\frac{(-1)^k}{k!} \Gamma(k+\beta-n+1) (1+\nu-n)_k \\
& \hspace{74mm} \times(m-2-k)! (1-z)^k \\
&\hspace{12mm} + \frac{(-1)^{m}}{\Gamma(\nu-n+1)} (1-z)^{m-1} \sum_{k=0}^{\infty} \frac{(\sigma)_k (\beta+\sigma-\nu)_k}{(k+\sigma-\nu+n-1)! k!} \\
&\hspace{44mm}\times \left[\log(1-z)-\psi(m+k)+\psi(m-n+\nu+k)\right.\\
&\hspace{54mm} \left. + \psi(k+\beta+m-n)-\psi(k+1)\right] (1-z)^k\\
&\hspace{12mm} - \sum_{l=0}^{n-2}\frac{1}{l!} \left[\sum_{k=l}^{n-2}\frac{(\nu-n+1)_{k} \Gamma(k-n+\beta+1)}{\Gamma(k-2n+\beta+m+\nu+1) (k-l)!}\right] (z-1)^l
\end{split}
\end{equation}
for $|1-z|<1$, $|\mathrm{arg}(1-z)|<\pi$, $n\in\mathbb{Z}^+$, $\nu\notin\mathbb{Z}^+$, $m\in\mathbb{Z}^+$ and $\mathrm{Re}(\beta)>(n-m)$. For given fixed $n$ and $m$, the restrictions on the parameters may be lifted by analytic continuation.  The transformation equations \eqref{res2} and \eqref{res2x} appear to be new.

\subsubsection{Case $(n+\sigma-\nu)\in\mathbb{Z}_0^-$} Under this condition, the integrals converges for $k=0,\dots, (\nu-\sigma-n)$ and the corresponding finite part integrals are just the regular values. Then the infinite series in the first term splits into two terms, 
\begin{equation}
\begin{split}
&\sum_{k=0}^{\infty} \binom{-\beta-\sigma+\nu}{k} \bbint{0}{\infty} \frac{((1-x)^{-1}+s)^{n-\nu-1}}{s^{\sigma-\nu+n+k}}\,\mathrm{d}s\\
&\hspace{14mm}=\sum_{k=0}^{\nu-\sigma-n} \binom{-\beta-\sigma+\nu}{k} \int_{0}^{\infty} s^{(\nu-\sigma-n)-k} ((1-x)^{-1}+s)^{n-\nu-1}\,\mathrm{d}s\\
& \hspace{18mm}+ \sum_{k=0}^{\infty} \binom{-\beta-\sigma+\nu}{k+\nu-\sigma-n+1}\bbint{0}{\infty} \frac{((1-x)^{-1}+s)^{n-\nu-1}}{s^{k+1}}\,\mathrm{d}s 
\end{split}
\end{equation}
where the index of summation for the remaining terms involving finite part integrals have been shifted to arrive at the second term above. The integrals in the first term are given by 
\begin{equation}
\int_{0}^{\infty} s^{(\nu-\sigma-n)-k} ((1-x)^{-1}+s)^{n-\nu-1}\,\mathrm{d}s = \frac{\Gamma(1-k-n+\nu-\sigma) \Gamma(k+\sigma)}{\Gamma(1-n+\nu)} (1-x)^{k+\sigma},
\end{equation}
for $k=0,\dots,(\nu-\sigma-n)$, $\mathrm{Re}(\sigma)>0$ and $\mathrm{Re}(\nu-n)>0$, which we obtain from \eqref{finiteres} with the substitutions $\lambda\rightarrow(-\nu+\sigma+n+k)$, $\upsilon\rightarrow (-n+\nu+1)$, and $s\rightarrow (1-x)^{-1}$. The conditions $\mathrm{Re}(\sigma)>0$ and $\mathrm{Re}(\nu-n)>0$ are there to ensure integrability at infinity for all $k$ in the given range. On the other hand, the fintie part integrals in the second term are given by
\begin{equation}
\begin{split}
\bbint{0}{\infty} \frac{((1-x)^{-1}+s)^{n-\nu-1}}{s^{k+1}}\,\mathrm{d}s =& \frac{(-1)^k (\nu+1-n)_k}{k!} (1-x)^{k+\nu+1-n}\\ &\times \left[-\ln(1-x)+\psi(k+1-\psi(k+\nu+1-n))\right]
\end{split}
\end{equation}
for $k\in\mathbb{Z}^+_0$, following from \eqref{FPIpole}. Integrability at infinity is guaranteed with the condition $\mathrm{Re}(\nu-n)>0$ for all $k$ in the range. Then  we have
\begin{equation}\label{bel1}
\begin{split}
&\sum_{k=0}^{\infty} \binom{-\beta-\sigma+\nu}{k} \bbint{0}{\infty} \frac{((1-x)^{-1}+s)^{n-\nu-1}}{s^{\sigma-\nu+n+k}}\,\mathrm{d}s\\
&\hspace{14mm}=\frac{1}{\Gamma(1-n+\nu)}\sum_{k=0}^{\nu-\sigma-n} \frac{(-1)^k}{k!} (\beta+\sigma-\nu)_k \Gamma(k+\sigma)\\ &\hspace{50mm}\times \Gamma(1+\nu-\sigma-n-k) (1-x)^{k+\sigma} \\
& \hspace{16mm}+ (-1)^{\nu-\sigma-n+1} \frac{\Gamma(\beta-n+1)}{\Gamma(\beta+\sigma-\nu)} \sum_{k=0}^{\infty} \frac{(\beta-n+1)_k (\nu+1-n)_k}{(k+\nu-\sigma-n+1)! k!} \\
&\hspace{30mm}\times\left[-\ln(1-x) + \psi(k+1)-\psi(k+\nu+1-n)\right] (1-x)^{k+\nu+1-n},
\end{split}
\end{equation}
which is valid for $0<x<1$, $n\in\mathbb{Z}^+$, $\mathrm{Re}(\sigma)>0$ and $\mathrm{Re}(\nu-n)>0$. 

Now the progenic finite part integral is still given by equation \eqref{mofpi}. However, since $(n+\sigma-\nu)$ is negative, $(2-n-\sigma+\nu+k)$ is positive for all $k=0,1,2,\dots$, so that the infinite series in the right hand side of \eqref{mofpi} does not splinter in two parts as in the previous case. Then the finite part integral in the series is given by \eqref{fpi2} for all $k$. Then
\begin{equation}\label{mofpix}
\begin{split}
	&\bbint{0}{1} \frac{s^{\beta-n} ((1-x)^{-1}-s)^{n-\nu-1}}{(1-s)^{\beta+\sigma-\nu}}\ln(s)\,\mathrm{d}s = (1-x)^{-(n-\nu-1)} \\
	&\hspace{14mm}\times (-1)^{\sigma-\nu} \frac{\pi \Gamma(\beta-n+1)}{\sin(\pi\beta) \Gamma(\beta+\sigma-\nu)} \sum_{k=0}^{\infty} \frac{(\nu+1-n)_k (\beta-n+1)_k}{(1-n-\sigma+\nu+k)! k!}\\
	&\hspace{14mm} \times\left[\psi(\beta-n+k+1)-\psi(\nu-\sigma-n+2+k)\right] (1-x)^k 
	\end{split}
\end{equation}
valid for $0<x<1$, $\mathrm{Re}(\beta)>(n-1)$, $\mathrm{\beta+\sigma-\nu}>0$ and $\beta\notin\mathbb{Z}$.  Substituting equations \eqref{bel1} and \eqref{mofpi} back into \eqref{fullint}, we obtain
\begin{equation}\label{keyxxx}
\begin{split}
&\int_0^{\infty} \frac{s^{\beta-n} ((1-x)^{-1}+ s)^{n-\nu-1}}{(s+1)^{\beta+\alpha-\nu}}\,\mathrm{d}s= \frac{1}{\Gamma(1-n+\nu)} \sum_{k=0}^{\nu-\sigma-n} \frac{(-1)^k}{k!} (\beta+\sigma-\nu)_k \\
&\hspace{24mm}\times\Gamma(k+\sigma) \Gamma(1+\nu-\sigma-n-k) (1-x)^{k+\sigma}\\
&\hspace{14mm}+(-1)^{\nu-\sigma-n} \frac{\Gamma(\beta-n+1)}{\Gamma(\beta+\sigma-\nu)} \sum_{k=0}^{\infty} \frac{(\beta-n+1)_k (\nu+1-n)_k}{(k+\nu-\sigma-n+1)! k!} \\
&\hspace{24mm}\times  \left[\ln(1-x)-\psi(k+1)+\psi(k+\nu+1-n)\right.\\
&\hspace{34mm}\left.+\psi(\beta-n+k+1)-\psi(\nu-\sigma-n+2+k)\right] (1-x)^{k-n+\nu+1}
\end{split}
\end{equation}
valid under the joint conditions for \eqref{bel1} and \eqref{mofpix}.

We substitute \eqref{keyxxx}, together with \eqref{leadint}, back into \eqref{3f2x}. Extending $x$ to complex values $z$ in the result, we obtain the transformation equation 
\begin{equation}
\begin{split}\label{bebebe}
&\frac{\Gamma(\beta) (\nu-n)_n  }{\Gamma(\beta+\sigma)(\nu-n)(n-1)!}z^{n-1}\, \pFq{3}{2}{\beta,\nu,1}{\beta+\sigma,n}{z}= \frac{1}{\Gamma(1-n+\nu)}(1-z)^{n-\nu-1+\sigma}\\
&\hspace{18mm} \times \sum_{k=0}^{\nu-\sigma-n} \frac{(-1)^k}{k!} (\beta+\sigma-\nu)_k (\sigma)_{k} \Gamma(1+\nu-\sigma-n-k) (1-z)^k\\
&\hspace{16mm} + (-1)^{\nu-\sigma-n} \frac{\Gamma(\beta-n+1)}{\Gamma(\beta+\sigma-\nu)\Gamma(\sigma)} \sum_{k=0}^{\infty} \frac{(\beta-n+1)_k (\nu+1-n)_k}{(k+\nu-\sigma-n+1)! k!}\\
&\hspace{18mm}\times \left[\log(1-z)-\psi(k+1)+\psi(k+\nu+1-n) + \psi(\beta-n+k+1)\right.\\
&\hspace{34mm}\left. -\psi(\nu-\sigma-n+2+k)\right] (1-z)^k\\
&\hspace{16mm} - \sum_{l=0}^{n-2}\frac{1}{l!} \left[\sum_{k=l}^{n-2}\frac{(\nu-n+1)_{k} \Gamma(k-n+\beta+1) k!}{\Gamma(k-n+\beta+\sigma+1) (k-l)!}\right] (z-1)^l  
\end{split}
\end{equation}
for $|1-z|<1$, $|\mathrm{arg}(1-z)|<\pi$, $n\in\mathbb{Z}^+$, $\nu\notin\mathbb{Z}$, $(\nu-\sigma-n)\in\mathbb{Z}^+_0$, $\mathrm{Re}(\sigma)>0$, $\mathrm{Re}(\nu-n)>0$, $\mathrm{Re}(\beta)>(n-1)$, $\mathrm{Re}(\beta+\sigma-\nu)>0$ and $\beta\notin\mathbb{Z}$, where we have performed again an expansion about $z=1$ to obtain the last term. 

We specialize with $(\nu+\sigma-n)=-m$ for $m\in\mathbb{Z}^+$ and eliminate $\sigma$ in favor of $n$, $m$ and $\nu$. Introducing the substitution $\sigma=\nu-m-n$ back into \eqref{bebebe}, we obtain 
\begin{equation}
\begin{split}\label{bebebex}
&\frac{\Gamma(\beta) (\nu-n)_n  }{\Gamma(\beta+\nu-m-n)(\nu-n)(n-1)!}z^{n-1}\, \pFq{3}{2}{\beta,\nu,1}{\beta+\nu-m-n,n,n}{z}\\
&= \frac{(1-z)^{-m-1}}{\Gamma(1-n+\nu)}
 \sum_{k=0}^{m} \frac{(-1)^k}{k!} (\beta-m-n)_k (\nu-m-n)_{k} \Gamma(1+m-k) (1-z)^k\\
&\hspace{16mm} + (-1)^{m} \frac{\Gamma(\beta-n+1)}{\Gamma(\beta-m-n) \Gamma(\nu-m-n)} \sum_{K=0}^{\infty} \frac{(\beta-n+1)_k (\nu+1-n)_k}{(k+m+1)! k!}\\
&\hspace{18mm}\times \left[\log(1-z)-\psi(k+1)+\psi(k+\nu+1-n) + \psi(\beta-n+k+1)\right.\\
&\hspace{34mm}\left. -\psi(k+m+s)\right] (1-z)^k\\
&\hspace{16mm} - \sum_{l=0}^{n-2}\frac{1}{l!} \left[\sum_{k=l}^{n-2}\frac{(\nu-n+1)_{k} \Gamma(k-n+\beta+1)k!}{\Gamma(k-2n-m+\nu+\beta+1) (k-l)!}\right] (z-1)^l ,
\end{split}
\end{equation}
for $|1-z|<1$, $|\mathrm{arg}(1-z)|<\pi$, $n\in\mathbb{Z}^+$, $m\in\mathbb{Z}^+_0$, $\mathrm{Re}(\nu)>(m+n)$ and $\mathrm{Re}(\beta)>(m+n)$. For a given fixed $m$ and $n$, the restrictions of the parameters $\nu$ and $\beta$ may be lifted by analytic continuation provided $\nu\notin\mathbb{Z}$. The transformation equations \eqref{bebebe} and \eqref{bebebex} appear to be new. 

\section{Conclusion}\label{conclusion}
We have elaborated and extended the method of finite part integration in the presence of competing singularities in the complex plane, and have applied the method to Stieltjes integral representations of the Gauss function and the generalized hypergeometric function $_3F_2$ that involve Stieltjes transforms of functions with complex extensions that have singularities in the complex plane. Finite-part integration of a Stieltjes integral representation of the Gauss hypergeometric function $_2F_1$ has led to transformation equations that reproduce known transformation equations of the Gauss function. Also new transformation equations involving the generalized hypergeometric function $_3F_2$ were obtained by direct finite part integration of an integral representation of $_3F_2$. The work here does not only extend the method of finite part integration but it also sets the framework upon which further properties of the generalized hypergeometric functions can be investigated through their generalized Stieltjes integral representations \cite{saxena,karp2}. In fact, a large class of transformation equations for the generalized hypergeometric function can already be derived by direct finite part integration of the integral representation given by \eqref{general}. The only hurdle to get through in performing finite part integration is the evaluation of the finite part integrals that arise. Currently our method of extracting finite parts has been limited to the application of the canonical definition of the finite part integral which may be intractable to implement under certain conditions. This calls for the development of new methods of evaluating the finite part of divergent integrals without explicit reference to the canonical definition of the finite part integral. As in regular integration where there are various methods of evaluating a convergent integral, we expect that different complementary methods may be developed to separate the finite part of divergent integrals.

\section*{Appendix} 

We derive the identity \eqref{keykey} using the canonical infinite series representation of the Gauss hypergeometric function written as
\begin{equation}
\pFq{2}{1}{\nu, 1}{n}{z} = \sum_{k=0}^{\infty} \frac{(\nu)_{k} (1)_{k}}{(n)_{k}} \frac{z^k}{k!}.
\end{equation}
Using the fact that $(1)_k = k!$ and the following identities
\begin{equation}
(n)_k = \frac{\Gamma(k+n)}{\Gamma(n)}, \;\;\; (\nu)_k = \frac{(\nu-n)_{k+n}}{(\nu-n)_n},
\end{equation}
we will obtain 
\begin{equation}
\pFq{2}{1}{\nu, 1}{n}{z} = \frac{\Gamma(n)}{z^{n} (\nu-n)_n } \sum_{k=0}^{\infty} \frac{(\nu-n)_{k+n}}{\Gamma(k+n)} z^{k+n}.
\end{equation}
Shifting the interval of summation, we will have
\begin{equation} \label{6.3}
\pFq{2}{1}{\nu, 1}{n}{z} = \frac{\Gamma(n)}{z^{n} (\nu-n)_n } \sum_{k=n}^{\infty} \frac{(\nu-n)_{k}}{\Gamma(k)} z^{k}.
\end{equation}
Take note that 
\begin{equation} \label{6.4}
\sum_{k=1}^{\infty} \frac{(\nu-n)_k z^k}{\Gamma(k)} = (\nu-n) (1-z)^{n-\nu-1} z.
\end{equation}
Extracting right-hand side of \eqref{6.4} from \eqref{6.3} will yield to
\begin{equation} \label{6.6}
\pFq{2}{1}{\nu, 1}{n}{z} = \frac{\Gamma(n)}{z^{n} (\nu-n)_n } \left( \sum_{k=1}^{\infty} \frac{(\nu-n)_{k}}{\Gamma(k)} z^{k} - \sum_{k=1}^{n-1} \frac{(\nu-n)_{k}}{\Gamma(k)} z^{k} \right).
\end{equation}
Substituting \eqref{6.4} to \eqref{6.6}, we obtain
\begin{equation} 
\pFq{2}{1}{\nu, 1}{n}{z} = \frac{\Gamma(n)}{z^{n} (\nu-n)_n } \left( (\nu-n) (1-z)^{n-\nu-1} z - \sum_{k=1}^{n-1} \frac{(\nu-n)_{k}}{\Gamma(k)} z^{k} \right).
\end{equation}
Making a shift on the summation interval and factoring out some common factor will result to
\begin{equation} 
\pFq{2}{1}{\nu, 1}{n}{z} = \frac{\Gamma(n)  (\nu-n)}{z^{n-1} (\nu-n)_n } \left( (1-z)^{n-\nu-1} - \sum_{k=0}^{n-2} \frac{(\nu-n+1)_{k}}{\Gamma(k+1)} z^{k} \right).
\end{equation}
The above equation is just the identity \eqref{keykey}.

\section*{Acknowledgement} 
The authors acknowledge the Office of the Chancellor of the University of the Philippines Diliman, through the Office of the Vice Chancellor for Research and Development, for funding support through the Outright Research Grant Project No. 191915.

\end{document}